\documentclass[11pt,a4paper]{article}

\usepackage[margin=20mm,top=22mm,bottom=24mm]{geometry}
\usepackage[T1]{fontenc}
\usepackage[utf8]{inputenc}
\usepackage{lmodern}
\usepackage{microtype}
\usepackage{amsmath,amssymb,amsthm,mathtools}
\usepackage{graphicx}
\usepackage[authoryear,square]{natbib}
\usepackage{placeins}
\usepackage{tikz}
\usetikzlibrary{arrows.meta,positioning,calc}
\usepackage[hidelinks]{hyperref}
\AddToHook{env/thebibliography/begin}{\interlinepenalty=10000\raggedbottom}
\hypersetup{
  pdftitle={Lyapunov--Schmidt reduction for periodic orbits of implicit discrete dynamical systems},
  pdfauthor={Henrique M. Oliveira}
}

\theoremstyle{plain}
\newtheorem{theorem}{Theorem}[section]
\newtheorem{proposition}[theorem]{Proposition}

\newtheorem{corollary}[theorem]{Corollary}
\theoremstyle{definition}
\newtheorem{definition}[theorem]{Definition}
\newtheorem{assumption}[theorem]{Assumption}
\theoremstyle{remark}
\newtheorem{remark}[theorem]{Remark}

\newcommand{\R}{\mathbb{R}}
\newcommand{\C}{\mathbb{C}}

\newcommand{\Id}{\mathrm{Id}}
\newcommand{\rank}{\operatorname{rank}}
\newcommand{\Range}{\operatorname{Range}}
\newcommand{\Ker}{\operatorname{Ker}}
\newcommand{\spec}{\operatorname{spec}}
\newcommand{\Sone}{\mathbb{S}^1}
\newcommand{\diag}{\operatorname{diag}}
\newcommand{\transp}{^{\mathsf T}}

\begin{document}

\title{Lyapunov--Schmidt reduction for periodic orbits of implicit discrete dynamical systems}

\author{Henrique M. Oliveira\\[0.5em]
\small Department of Mathematics, Instituto Superior T\'ecnico, Universidade de Lisboa\\
\small Centro de An\'alise Matem\'atica, Geometria e Sistemas Din\^amicos (CAMGSD)\\
\small Av. Rovisco Pais, 1049-001 Lisboa, Portugal\\[0.3em]
\small \href{mailto:henrique.m.oliveira@tecnico.ulisboa.pt}{henrique.m.oliveira@tecnico.ulisboa.pt}}
\date{}

\maketitle

\begin{abstract}
We derive bifurcation equations for periodic orbits of autonomous and
periodically non-autonomous implicit discrete systems in arbitrary
finite dimension. Cyclic Lyapunov--Schmidt reduction gives recursive
coefficients and nondegeneracy and transversality conditions at the
multipliers $+1$ and $-1$. We prove equivalence under regular changes of
implicit representation and exact invariance under cyclic relabelling
when the reduction data are transported consistently.
\end{abstract}



\section{Introduction}
\label{sec:introduction}

\subsection{Origins of the problem}
In 2016 and 2018 I published two papers on bifurcation equations
\citep{Oliveira2016,Oliveira2018}. The first established invariance under
cyclic reordering for periodically non-autonomous systems of explicit
maps, leaving the multiplier $-1$ untreated. The second derived
bifurcation equations and nondegeneracy and transversality conditions
for implicit discrete systems. Those formulae grew rapidly with the
order of degeneracy, as successive compositions required increasingly
long derivative expressions.

I returned to these questions by describing each periodic orbit through
all its states simultaneously. The cyclic equations treat autonomous
and periodically non-autonomous implicit systems in arbitrary finite
dimension; period doubling is included by writing the orbit twice.

\subsection{Motivation and relation with previous works}

For states $x_n\in\R^d$ and parameters $\lambda\in\R^m$, consider
\begin{equation}
    F_n(x_n,x_{n+1},\lambda)=0.
    \label{eq:implicit-map}
\end{equation}
The system is autonomous when $F_n$ is independent of $n$, and
periodically non-autonomous when $F_{n+p}=F_n$ for a positive integer
$p$. We study orbits closing after $q$ steps: $q$ is arbitrary in
the autonomous case and $q=Np$, with $N$ a positive integer, in the $p$-periodic non-autonomous
case. This prescribed length need not be the least period of the orbit.
We assume that the derivative with respect to the next state is
invertible along the orbit. The implicit function theorem then gives
locally unique forward steps. Their composition over $q$ steps has a
linearisation, the \emph{monodromy matrix}; its eigenvalues are the
\emph{multipliers}.

Implicit relations occur naturally in numerical integration. Backward
Euler, trapezoidal and implicit Runge--Kutta methods determine the next
state by solving nonlinear equations
\citep{HairerWanner1996,BrenanCampbellPetzold1996}. Implicit mappings
also arise in mechanical models, including machining processes
\citep{Luo2015,Luo2015PeriodicFlows,KissBachrathy2022}. In economics,
relations between present and anticipated future states may admit
several forward branches. Their global dynamics has been studied
through inverse limits by Medio and Raines \citeyearpar{MedioRaines2007},
Kennedy, Stockman and Yorke
\citeyearpar{KennedyStockmanYorke2007,KennedyStockmanYorke2008}, and Kennedy
and Stockman \citeyearpar{KennedyStockman2008}. The present analysis concerns
local bifurcations along regular periodic orbits.

Normal forms and bifurcations of difference-algebraic equations
were studied by Beardmore and Webster \citeyearpar{BeardmoreWebster2008},
while de la Llave and Lomel\'i \citeyearpar{DeLaLlaveLomeli2012}
developed a parameterisation method for invariant manifolds of
difference equations. Multiple-shooting methods for invariant
manifolds attached to periodic orbits were developed for maps by
Gonzalez and Mireles James \citeyearpar{GonzalezMirelesJames2017} and for
implicitly defined systems by Timsina and Mireles James
\citeyearpar{TimsinaMirelesJames2022}. Here the emphasis is instead on
bifurcation equations obtained by reducing the cyclic relations
defining the periodic orbit.

For explicit maps, the local theory of bifurcations and normal forms
is well established
\citep{ArrowsmithPlace1990,Elaydi2005,KuznetsovBook,KuznetsovMeijer2005,GovaertsEtAl2007,BeynEtAl2002}. For a systematic treatment of higher-codimension bifurcations
of maps, see also Kuznetsov and Meijer
\citeyearpar{KuznetsovMeijer2019}.
Lyapunov--Schmidt reduction of periodic-point equations has been used
to retain the symmetries of periodic sequences by Ciocci
\citeyearpar{Ciocci2004} and Vanderbauwhede \citeyearpar{Vanderbauwhede2010}; the
relevant singularity theory is developed in
\citet{GolubitskySchaeffer1985}, and the equivalence of nonlinear
reductions is analysed in \citet{JepsonSpence1989}. Kiss and Bachrathy
\citeyearpar{KissBachrathy2022} subsequently computed normal-form
coefficients for a scalar implicit model of milling. Their account of
previous work does not mention the arbitrary-period bifurcation
equations and nondegeneracy and transversality conditions derived in
\citet{Oliveira2018}. The two papers address different questions:
periodic-orbit bifurcation criteria in the earlier work and
normal-form coefficients for a particular engineering model in the
later one.

The cyclic construction gives a spectral identity relating the orbit
linearisation to the monodromy matrix, a recursion for reduced
coefficients of arbitrary order, and comparison theorems for the
reduced equations. A regular change of implicit representation
multiplies the scalar equation by a smooth nonvanishing factor.
Cyclic relabelling preserves it exactly when the kernel vectors are
transported with the orbit. The doubled equations extend both results
to ordinary period doubling and the codimension-two degenerate flip. These are normal forms
of bifurcation equations under the equivalence defined below.

The paper is organised as follows. Section~\ref{sec:cyclic} introduces
the cyclic orbit equations and their spectral linearisation.
Section~\ref{sec:representation} establishes independence of the
implicit representation. Section~\ref{sec:LS} constructs the
Lyapunov--Schmidt reduction and its recursive coefficients.
Section~\ref{sec:plus-one} treats the bifurcations at the multiplier
$+1$, and Section~\ref{sec:doubling} treats period doubling.
Section~\ref{sec:cyclic-invariance} establishes invariance under
cyclic relabelling, including for periodically non-autonomous systems.
Sections~\ref{sec:examples} and \ref{sec:discussion} give examples and conclusions.

\section{Periodic Orbits as Cyclic Implicit Systems}
\label{sec:cyclic}

Let $d,m$ be positive integers, let $\Lambda\subset\R^m$ be open, and let $u,v\in\R^d$
denote the current and next state. Write $I$ for an identity matrix
of the required size. For a smooth relation $F$, write
\begin{equation}
    F(u,v,\lambda)=0,\qquad \lambda\in\Lambda.
    \label{eq:two-state-relation}
\end{equation}
The matrices $D_1F$ and $D_2F$ are its derivatives with respect to
$u$ and $v$, respectively. The argument subscripts do not denote time
indices. Invertibility of $D_2F$ gives a locally unique forward step;
we use this convention throughout \citep{Luo2015,Luo2015PeriodicFlows}.

Fix a positive integer $q$, open sets $U_j\subset\R^d$ with
$U_{j+q}=U_j$, and smooth relations
$F_j:U_j\times U_{j+1}\times\Lambda\to\R^d$, extended by
$F_{j+q}=F_j$. Here and below smooth means $C^\infty$.
The orbit equations are
\begin{equation}
    F_j(x_j,x_{j+1},\lambda)=0,\qquad x_{j+q}=x_j,
    \quad j\in\mathbb Z.
    \label{eq:periodic-implicit-system}
\end{equation}
For an autonomous system, the $F_j$ are restrictions of one relation
and $q$ is arbitrary. For a periodically non-autonomous system with
coefficient period $p$, choose $q=Np$ with $N$ a positive integer and use the
$p$-periodic defining sequence in \eqref{eq:periodic-implicit-system}.
The coefficient period $p$ and the prescribed orbit length $q$ need not coincide.

A \emph{$q$-periodic solution} satisfies
\eqref{eq:periodic-implicit-system}; its \emph{least period} is the
smallest positive integer $r$ for which $x_{j+r}=x_j$ for every $j$.
We allow repeated orbits, so this least period can divide $q$.
Collect the states in
\[
    X=(x_0,\ldots,x_{q-1})\in
    U_0\times\cdots\times U_{q-1}\subset\R^{dq}.
\]
The periodic-orbit problem has $dq$ equations in $dq$ state variables.

\begin{definition}
For $q\geq2$, the \emph{cyclic orbit map} of
\eqref{eq:periodic-implicit-system} is the map
\begin{equation}
    \Phi(X,\lambda)=
    \begin{pmatrix}
        F_0(x_0,x_1,\lambda)\\
        F_1(x_1,x_2,\lambda)\\
        \vdots\\
        F_{q-2}(x_{q-2},x_{q-1},\lambda)\\
        F_{q-1}(x_{q-1},x_0,\lambda)
    \end{pmatrix}.
    \label{eq:cyclic-map}
\end{equation}
For $q=1$, the definition is simply
$\Phi(x_0,\lambda)=F_0(x_0,x_0,\lambda)$.
A vector $X$ represents a $q$-periodic orbit if and only if
$\Phi(X,\lambda)=0$.
\end{definition}

Fix a periodic solution $(X_*,\lambda_*)$ and write
$X_*=(x_0^*,\ldots,x_{q-1}^*)$, with $x_q^*=x_0^*$.
The derivatives of the individual relations at this orbit are
\begin{equation}
    A_j=D_1F_j(x_j^*,x_{j+1}^*,\lambda_*),
    \qquad
    B_j=D_2F_j(x_j^*,x_{j+1}^*,\lambda_*).
    \label{eq:AjBj}
\end{equation}
Each $A_j$ and $B_j$ is a real $d\times d$ matrix.
Our standing assumption is that the second-argument derivatives are
invertible:
\begin{equation}
    \det B_j\neq0,
    \qquad j=0,\ldots,q-1.
    \label{eq:regularity}
\end{equation}
This is the regularity condition for the chosen forward evolution.

Under \eqref{eq:regularity}, the implicit function theorem gives a
locally defined map $f_j$ satisfying
\[
    F_j(u,f_j(u,\lambda),\lambda)=0,
    \qquad f_j(x_j^*,\lambda_*)=x_{j+1}^*.
\]
Differentiating the first identity with respect to $u$ at
$(x_j^*,\lambda_*)$ gives
\[
    A_j+B_jD_uf_j(x_j^*,\lambda_*)=0.
\]
Thus the derivative of the forward step is
\begin{equation}
    J_j=D_uf_j(x_j^*,\lambda_*)=-B_j^{-1}A_j.
    \label{eq:step-derivative}
\end{equation}
The product of these derivatives over one full period is the
\emph{monodromy matrix}
\begin{equation}
    M=J_{q-1}J_{q-2}\cdots J_0.
    \label{eq:monodromy}
\end{equation}
It is the derivative at $x_0^*$ of the local return map
$f_{q-1}\circ\cdots\circ f_0$, with $\lambda=\lambda_*$.
Its spectrum is denoted by $\spec(M)$. Its eigenvalues, counted
with algebraic multiplicity, are the \emph{multipliers} of the periodic orbit \citep{DaCunhaDavis2011,Teschl2000}.
In particular, the multipliers $+1$ and $-1$ are the critical values
for the local bifurcations studied later.

To see how this works, let $v_j\in\C^d$ be a perturbation of $x_j^*$.
Working over $\C$ includes non-real multipliers. Linearising the
$j$th implicit equation gives
\[
    A_jv_j+B_jv_{j+1}=0,
    \qquad\text{or, equivalently,}\qquad
    v_{j+1}=J_jv_j.
\]
After $q$ steps, every such perturbation satisfies
$v_q=Mv_0$. To find the multipliers, we seek nonzero initial
perturbations $v_0$ that are multiplied by a scalar $\zeta$ after
one full period:
\begin{equation}
    v_q=\zeta v_0,
    \qquad \zeta\in\C.
    \label{eq:twisted-closure}
\end{equation}
Together with the linearised evolution, this is precisely the
eigenvalue condition $Mv_0=\zeta v_0$. A general perturbation need
not be periodic, even though the reference orbit satisfies $x_q^*=x_0^*$.
We call \eqref{eq:twisted-closure} a \emph{twisted} closing condition:
the initial perturbation is reproduced up to the factor $\zeta$
after one circuit. This is a Floquet-type condition for a periodic
difference system \citep{DaCunhaDavis2011,Teschl2000}. For $\zeta=1$
the perturbation is periodic; for $\zeta=-1$ it is anti-periodic.

\begin{definition}[Twisted cyclic operator]
\label{def:twisted}
For $\zeta\in\C$, let
$L_\zeta:(\C^d)^q\longrightarrow(\C^d)^q$
act on $V=(v_0,\ldots,v_{q-1})$ through the linearised orbit
equations, with the closing condition \eqref{eq:twisted-closure}:
\begin{equation}
    (L_\zeta V)_j=A_jv_j+B_jv_{j+1},
    \qquad j=0,\ldots,q-2,
    \label{eq:twisted-first}
\end{equation}
\begin{equation}
    (L_\zeta V)_{q-1}
       =A_{q-1}v_{q-1}+\zeta B_{q-1}v_0.
    \label{eq:twisted-last}
\end{equation}
For $q=1$, this reads
\begin{equation}
    L_\zeta=A_0+\zeta B_0.
    \label{eq:p1-twisted}
\end{equation}
In particular, $L_1=D_X\Phi(X_*,\lambda_*)$.
\end{definition}

For $q\geq2$, the block matrix of $L_\zeta$ is
\begin{equation}
L_\zeta=
\begin{pmatrix}
 A_0 & B_0 & 0 & \cdots & 0\\
 0 & A_1 & B_1 & \ddots & \vdots\\
 \vdots & \ddots & \ddots & \ddots & 0\\
 0 & \cdots & 0 & A_{q-2} & B_{q-2}\\
 \zeta B_{q-1} & 0 & \cdots & 0 & A_{q-1}
\end{pmatrix}.
\label{eq:twisted-matrix}
\end{equation}
The lower-left block records the closing factor $\zeta$. The
singularities of $L_\zeta$ determine the multipliers.

\begin{theorem}[Twisted cyclic--monodromy correspondence]
\label{thm:twisted-monodromy}
Under \eqref{eq:regularity}, for every $\zeta\in\C$ projection onto
the first component gives a linear isomorphism
\begin{equation}
    \Ker L_\zeta\longrightarrow\Ker(M-\zeta I).
    \label{eq:kernel-isomorphism}
\end{equation}
Moreover,
\begin{equation}
    \det L_\zeta
    =(-1)^{d(q-1)}
      \left(\prod_{j=0}^{q-1}\det B_j\right)
      \det(\zeta I-M).
    \label{eq:det-identity}
\end{equation}
\end{theorem}

In particular,
\begin{equation}
    \dim\Ker L_\zeta=\dim\Ker(M-\zeta I),
    \label{eq:kernel-dim}
\end{equation}
so $L_\zeta$ is singular exactly at a multiplier.

\begin{proof}
Suppose that $V=(v_0,\ldots,v_{q-1})\in\Ker L_\zeta$.
Since every $B_j$ is invertible, the first $q-1$ equations give
\[
    v_{j+1}=J_jv_j,
    \qquad j=0,\ldots,q-2.
\]
Hence $v_j=J_{j-1}\cdots J_0v_0$ for $1\leq j\leq q-1$.
The final equation is
\[
    B_{q-1}(\zeta v_0-J_{q-1}v_{q-1})=0,
\]
which is equivalent to $Mv_0=\zeta v_0$. Conversely, given any
$v_0\in\Ker(M-\zeta I)$, the recurrence determines all the other
$v_j$ and satisfies the final equation. This proves the isomorphism
and the equality of nullities. For $q=1$ the same conclusion follows
directly from $L_\zeta=B_0(\zeta I-J_0)$.

To establish the determinant identity, factor $B_j$ from block row
$j$. For $q\geq2$ this gives
\[
    \det L_\zeta
    =\left(\prod_{j=0}^{q-1}\det B_j\right)\det K_\zeta,
\]
where
\[
K_\zeta=
\begin{pmatrix}
 -J_0&I&0&\cdots&0\\
 0&-J_1&I&\ddots&\vdots\\
 \vdots&\ddots&\ddots&\ddots&0\\
 0&\cdots&0&-J_{q-2}&I\\
 \zeta I&0&\cdots&0&-J_{q-1}
\end{pmatrix}.
\]
Move the first block column to the last position. This permutation
contributes the sign $(-1)^{d(q-1)}$. The upper-left block of the
rearranged matrix is lower triangular with identity blocks on the
diagonal, and its Schur complement is
\[
    \zeta I-J_{q-1}\cdots J_0=\zeta I-M.
\]
Thus $\det K_\zeta=(-1)^{d(q-1)}\det(\zeta I-M)$,
as required. For $q=1$ the formula follows from
$L_\zeta=B_0(\zeta I-J_0)$.
\end{proof}

\begin{corollary}
\label{cor:critical-multipliers}
Under \eqref{eq:regularity},
\[
    \Ker D_X\Phi(X_*,\lambda_*)\neq\{0\}
    \quad\Longleftrightarrow\quad
    1\in\spec(M),
\]
and
\[
    \Ker L_{-1}\neq\{0\}
    \quad\Longleftrightarrow\quad
    -1\in\spec(M).
\]
\end{corollary}

\begin{remark}
When $d=1$ and $F_j=F$, the derivative of the local return map,
formed by following the branches $f_j$ along the orbit, is
\[
    (f_{q-1}\circ\cdots\circ f_0)'(x_0^*)
    =(-1)^q\prod_{j=0}^{q-1}
      \frac{D_1F(x_j^*,x_{j+1}^*,\lambda_*)}
           {D_2F(x_j^*,x_{j+1}^*,\lambda_*)}.
\]
This is the scalar multiplier formula of \citet{Oliveira2018};
Theorem~\ref{thm:twisted-monodromy} gives its vector-valued form.
\end{remark}

\section{Invariance under Regular Implicit Representations}
\label{sec:representation}

In the scalar theory of \citet{Oliveira2018}, bifurcation conditions
are calculated from derivatives of the defining relation. Different
regular relations can describe the same local dynamics while having
different derivatives. We must therefore establish that the resulting
bifurcation criteria are independent of this choice. We first express
one defining relation as a smooth invertible matrix times the other
and prove that their linearised dynamics agree. Section~\ref{sec:LS}
uses this factorisation to show that the scalar reduced equations
differ by a smooth nonvanishing factor. Their degeneracy and
transversality conditions are consequently unchanged, although the
numerical values of their coefficients may differ.

Consider two smooth functions $F,\widetilde F$ from a neighbourhood of
$(u_*,v_*,\lambda_*)$ in $\R^d\times\R^d\times\Lambda$ to $\R^d$.
We call them \emph{equivalent local representations} if they vanish at
$(u_*,v_*,\lambda_*)$ and their zero sets agree in a neighbourhood of
that point. Both representations are assumed to be forward regular, in
the sense of Section~\ref{sec:cyclic}: $D_2F$ and $D_2\widetilde F$
are invertible at the solution.

Since $D_2F(u_*,v_*,\lambda_*)$ is invertible, the map
\[
    (u,v,\lambda)\longmapsto(u,F(u,v,\lambda),\lambda)
\]
is a local diffeomorphism. Write its inverse as
$(u,z,\lambda)\mapsto(u,\vartheta(u,z,\lambda),\lambda)$ and define
\[
\begin{aligned}
    H(u,z,\lambda)&=\widetilde F(u,\vartheta(u,z,\lambda),\lambda),\\
    \Gamma(u,v,\lambda)
      &=\int_0^1 D_zH(u,tF(u,v,\lambda),\lambda)\,dt.
\end{aligned}
\]

\begin{theorem}[Factorisation of regular representations]
\label{thm:representation-factorization}
The matrix-valued function $\Gamma$ is smooth and invertible near
$(u_*,v_*,\lambda_*)$, and
\begin{equation}
    \widetilde F(u,v,\lambda)=\Gamma(u,v,\lambda)F(u,v,\lambda).
    \label{eq:representation-factorization}
\end{equation}
\end{theorem}

\begin{proof}
The common zero set gives $H(u,0,\lambda)=0$. The fundamental theorem
of calculus therefore gives
\[
    H(u,z,\lambda)
      =\left(\int_0^1D_zH(u,tz,\lambda)\,dt\right)z.
\]
Substitution of $z=F(u,v,\lambda)$ proves
\eqref{eq:representation-factorization}. Smoothness follows from the
definition, and
\[
    \Gamma(u_*,v_*,\lambda_*)
      =D_2\widetilde F(u_*,v_*,\lambda_*)
        [D_2F(u_*,v_*,\lambda_*)]^{-1}
\]
is invertible. It remains invertible in a neighbourhood.
\end{proof}

Apply the factorisation to each relation along the $q$-periodic orbit
from Section~\ref{sec:cyclic}, at the corresponding points of the orbit
$X_*=(x_0^*,\ldots,x_{q-1}^*)$. Let $\Gamma_j$ be the factor for the
$j$th relation and define
\begin{equation}
    \mathcal T(X,\lambda)
    =\diag\bigl(\Gamma_0(x_0,x_1,\lambda),\ldots,
                  \Gamma_{q-1}(x_{q-1},x_0,\lambda)\bigr),
    \label{eq:block-representation}
\end{equation}
where $X=(x_0,\ldots,x_{q-1})$. The block-diagonal matrix
$\mathcal T$ acts on the $dq$-dimensional orbit space.
Thus the two cyclic orbit maps satisfy the identity
\begin{equation}
    \widetilde\Phi(X,\lambda)
    =\mathcal T(X,\lambda)\Phi(X,\lambda).
    \label{eq:cyclic-representation}
\end{equation}
In particular, they have the same local periodic-orbit solutions for
every nearby parameter value. Write
$\mathcal T_*=\mathcal T(X_*,\lambda_*)$.
Tildes denote the matrices and operators constructed from
$\widetilde F_j$ by the same definitions.

\begin{theorem}[Invariance of the linearised dynamics]
\label{thm:representation-invariance}
Suppose that $F_j$ and $\widetilde F_j$ are equivalent forward-regular
representations along a $q$-periodic orbit. For every $\zeta\in\C$,
\begin{equation}
    \widetilde L_\zeta=\mathcal T_*L_\zeta,
    \qquad
    \Ker\widetilde L_\zeta=\Ker L_\zeta.
    \label{eq:L-representation}
\end{equation}
Moreover,
\begin{equation}
    \widetilde J_j=J_j\quad (j=0,\ldots,q-1),
    \qquad
    \widetilde M=M.
    \label{eq:CM-representation}
\end{equation}
The multipliers and their eigenspaces are therefore independent of the
regular implicit representation.
\end{theorem}

\begin{proof}
Write $\Gamma_j^*=\Gamma_j(x_j^*,x_{j+1}^*,\lambda_*)$.
At this solution, differentiation of $\widetilde F_j=\Gamma_jF_j$
with respect to either state argument gives
\[
    \widetilde A_j=\Gamma_j^*A_j,
    \qquad
    \widetilde B_j=\Gamma_j^*B_j.
\]
Hence the forward-step derivatives satisfy
\[
    \widetilde J_j
       =-(\widetilde B_j)^{-1}\widetilde A_j
       =-B_j^{-1}A_j=J_j.
\]
Their products give $\widetilde M=M$. Multiplying each block row of
$L_\zeta$ by $\Gamma_j^*$ gives the first identity in
\eqref{eq:L-representation}; its kernel identity follows because
$\mathcal T_*$ is invertible.
\end{proof}

\section{Lyapunov--Schmidt Reduction at a Simple Multiplier \texorpdfstring{$+1$}{+1}}
\label{sec:LS}

Let $(X_*,\lambda_*)$ be a periodic solution whose monodromy matrix
has an algebraically simple multiplier $+1$. Lyapunov--Schmidt
reduction separates orbit displacements along the kernel of the
linearisation from those in a complementary space. The complementary
equations are solved by the implicit function theorem, leaving one
scalar equation for the nearby periodic orbits. We carry out the
reduction in the full orbit space $\R^{dq}$.

Use the real restriction of $L_1$ and put
\[
    L=L_1=D_X\Phi(X_*,\lambda_*).
\]
By Theorem~\ref{thm:twisted-monodromy}, $\Ker L$ is one-dimensional.
Since $L$ is a square $dq\times dq$ matrix, the same holds for its
left kernel:
\begin{equation}
    \dim\Ker L=\dim\Ker L\transp=1.
    \label{eq:simple-kernel}
\end{equation}
Use the Euclidean inner product $\langle\cdot,\cdot\rangle$ and norm
$\|\cdot\|$ on $\R^{dq}$. Choose unit vectors
\begin{equation}
    V\in\Ker L,\qquad W\in\Ker L\transp,
    \qquad \|V\|=\|W\|=1.
    \label{eq:left-right-vectors}
\end{equation}
The vector $V$ spans the critical direction in the orbit
space, while $W$ is orthogonal to the range of $L$. In particular,
\[
    \Range L=W^\perp,
    \qquad
    V^\perp=\{Z\in\R^{dq}:\langle V,Z\rangle=0\}.
\]
The restriction $L|_{V^\perp}:V^\perp\to W^\perp$ is an isomorphism:
it is injective because $\Ker L=\operatorname{span}\{V\}$, and both
spaces have dimension $dq-1$. Define the projection onto the range
and the corresponding reduced inverse by
\begin{equation}
    Q=I-WW\transp,\qquad
    \mathcal R=(L|_{V^\perp})^{-1}Q:
       \R^{dq}\longrightarrow V^\perp.
    \label{eq:reduced-inverse}
\end{equation}
Thus $\mathcal R b$ is the unique vector in $V^\perp$ whose image
under $L$ is $Qb$.

Every $X$ near $X_*$ has a unique decomposition
\begin{equation}
    X=X_*+sV+Z,\qquad s\in\R,\quad Z\in V^\perp,
    \label{eq:LS-splitting}
\end{equation}
where $s=\langle V,X-X_*\rangle$. We first solve the $dq-1$
components of the cyclic equation belonging to $W^\perp$:
\begin{equation}
    Q\Phi(X_*+sV+Z,\lambda)=0.
    \label{eq:range-equation}
\end{equation}
Its derivative with respect to $Z$ at
$(s,Z,\lambda)=(0,0,\lambda_*)$ is $L|_{V^\perp}$, which is invertible.
The implicit function theorem therefore gives a unique smooth map
$Z=Z(s,\lambda)\in V^\perp$, with $Z(0,\lambda_*)=0$. Set
\[
    Y(s,\lambda)=X_*+sV+Z(s,\lambda).
\]
The range equation implies that $\Phi(Y(s,\lambda),\lambda)$ is a
multiple of $W$. Its scalar component is
\begin{equation}
    \psi(s,\lambda)
    =\langle W,\Phi(Y(s,\lambda),\lambda)\rangle.
    \label{eq:reduced-scalar}
\end{equation}
In particular,
\[
    \Phi(Y(s,\lambda),\lambda)=\psi(s,\lambda)W.
\]
Differentiation of \eqref{eq:range-equation} with respect to $s$ at
the base point gives $LZ_s(0,\lambda_*)=0$. Since $Z_s$ belongs to
$V^\perp$, we obtain
\begin{equation}
    Z_s(0,\lambda_*)=0,\qquad
    \psi(0,\lambda_*)=\psi_s(0,\lambda_*)=0.
    \label{eq:LS-base-derivatives}
\end{equation}

\begin{theorem}[Cyclic Lyapunov--Schmidt reduction]
\label{thm:LS}
Under \eqref{eq:regularity} and \eqref{eq:simple-kernel}, the nearby
periodic-orbit solutions are exactly
\[
    X=Y(s,\lambda),\qquad \psi(s,\lambda)=0.
\]
\end{theorem}

\begin{proof}
The decomposition \eqref{eq:LS-splitting} assigns a unique pair
$(s,Z)$ to every nearby $X$. By the implicit function theorem, the
range equation holds exactly when $Z=Z(s,\lambda)$. On that graph,
$\Phi=\psi W$, and hence $\Phi=0$ exactly when $\psi=0$.
\end{proof}

Let $f_{j,\lambda}(u)=f_j(u,\lambda)$ and write
$P_\lambda=f_{q-1,\lambda}\circ\cdots\circ f_{0,\lambda}$
for the local return map. For $q\geq2$, define the residual vectors
\[
    e_j=x_{j+1}-f_j(x_j,\lambda)\in\R^d,
    \qquad 0\leq j<q-1,
\]
and set $e=(e_0,\ldots,e_{q-2})\in(\R^d)^{q-1}$.
Thus $e=0$ expresses the first $q-1$ evolution equations.
For $q=1$, the block $e$ is absent.

\begin{proposition}[Equivalence with the return equation]
\label{prop:return-equation}
Under \eqref{eq:regularity}, the cyclic orbit equations are locally
equivalent, by smooth changes of state variables and invertible
recombinations of equations, to
\[
    e=0,\qquad x_0-P_\lambda(x_0)=0.
\]
At a simple multiplier $+1$, the corresponding scalar
Lyapunov--Schmidt reductions differ only by a smooth local change
of critical coordinate and a smooth nonvanishing factor.
\end{proposition}

\begin{proof}
For $q\geq2$, the map $X\mapsto(x_0,e)$ is a local
diffeomorphism.
The fundamental theorem of calculus gives
\[
 F_j(u,v,\lambda)=C_j(u,v,\lambda)
              \bigl(v-f_j(u,\lambda)\bigr),
\]
where $C_j$ is smooth and invertible near the orbit. Therefore
an invertible recombination of the cyclic equations gives
\[
 e_0=\cdots=e_{q-2}=0,\qquad
 x_0-f_{q-1}(x_{q-1},\lambda)=0.
\]
In the new coordinates,
$f_{q-1}(x_{q-1},\lambda)-P_\lambda(x_0)$ vanishes when
$e=0$ and hence equals $H(x_0,e,\lambda)e$ for some smooth
matrix $H$. Adding $He$ to the last equation gives the claimed
system by an invertible block-triangular recombination. For $q=1$
the same conclusion follows directly from the factorisation of
$F_0$. Eliminating the regular variables $e$ leaves precisely
the fixed-point equation. The equivalence of the associated scalar
Lyapunov--Schmidt reductions follows from invariance under smooth
changes of variables and invertible recombinations of equations
\citep{JepsonSpence1989}.
\end{proof}

The comparison below concerns changes of implicit representation.
The same reduction will be used in Section~\ref{sec:cyclic-invariance}
to extend the cyclic invariance results established for explicit systems
\citep{Oliveira2016}.

Consider an equivalent forward-regular representation
$\widetilde\Phi=\mathcal T\Phi$ as in
\eqref{eq:cyclic-representation}. Use the same critical coordinate
$s$ and right-kernel vector $V$ for both reductions. For the second
representation choose the unit left-kernel vector
\begin{equation}
    \widetilde W
    =\frac{\mathcal T_*^{-\mathsf T}W}
           {\|\mathcal T_*^{-\mathsf T}W\|}.
    \label{eq:transformed-left-vector}
\end{equation}
Put $\widetilde Q=I-\widetilde W\widetilde W\transp$ and let
$\widetilde Z(s,\lambda)\in V^\perp$ be the unique smooth solution of
\[
    \widetilde Q\widetilde\Phi
       (X_*+sV+\widetilde Z(s,\lambda),\lambda)=0.
\]
The corresponding reduced equation is
\[
    \widetilde\psi(s,\lambda)
       =\bigl\langle\widetilde W,
          \widetilde\Phi(X_*+sV+\widetilde Z(s,\lambda),\lambda)
        \bigr\rangle.
\]

For $\Delta\in V^\perp$, define
\[
    \mathcal K(s,\lambda,\Delta)
      =\left.\int_0^1D_X\Phi(Y(s,\lambda)+\tau\Delta,\lambda)
        \,d\tau\right|_{V^\perp}.
\]
With $Y=Y(s,\lambda)$ and an auxiliary scalar $t$, consider
\begin{equation}
    \widetilde Q\,\mathcal T(Y+\Delta,\lambda)
      \bigl(tW+\mathcal K(s,\lambda,\Delta)\Delta\bigr)=0.
    \label{eq:comparison-range}
\end{equation}
At $(s,\lambda,t,\Delta)=(0,\lambda_*,0,0)$, the derivative with
respect to $\Delta$ is the isomorphism
$\widetilde L_1|_{V^\perp}:V^\perp\to\widetilde W^\perp$.
The implicit function theorem gives a unique smooth solution
$\Delta(s,\lambda,t)$, with $\Delta(s,\lambda,0)=0$.
Set
\[
    \Theta(s,\lambda,t)
      =\int_0^1\partial_t\Delta(s,\lambda,\tau t)\,d\tau,
    \qquad
    \Delta_\psi(s,\lambda)=\Delta(s,\lambda,\psi(s,\lambda)).
\]
Then $\Delta(s,\lambda,t)=t\Theta(s,\lambda,t)$. Define
\[
    c(s,\lambda)=
      \left\langle\widetilde W,
        \mathcal T(Y+\Delta_\psi,\lambda)
        \bigl(W+\mathcal K(s,\lambda,\Delta_\psi)
                    \Theta(s,\lambda,\psi(s,\lambda))\bigr)
      \right\rangle.
\]

\begin{theorem}[Invariance of the reduced equation]
\label{thm:reduced-representation}
Under the hypotheses of Theorem~\ref{thm:LS}, the function $c$ is
smooth and satisfies
\begin{equation}
    \widetilde\psi(s,\lambda)=c(s,\lambda)\psi(s,\lambda),
    \qquad c(0,\lambda_*)\ne0.
    \label{eq:scalar-representation}
\end{equation}
\end{theorem}

\begin{proof}
The fundamental theorem of calculus gives
\[
    \Phi(Y+\Delta,\lambda)
      =\psi(s,\lambda)W+\mathcal K(s,\lambda,\Delta)\Delta.
\]
Taking $t=\psi(s,\lambda)$ in \eqref{eq:comparison-range} therefore
gives the range equation for $\widetilde\Phi$ at $Y+\Delta_\psi$.
Uniqueness implies
\[
    \widetilde Z(s,\lambda)=Z(s,\lambda)+\Delta_\psi(s,\lambda).
\]
Substituting $\Delta_\psi=\psi\Theta(s,\lambda,\psi)$ into the scalar
component gives \eqref{eq:scalar-representation}. Smoothness follows
from the construction. Since $\widetilde W\transp\mathcal T_*L=0$,
\[
    c(0,\lambda_*)
      =\langle\widetilde W,\mathcal T_*W\rangle
      =\|\mathcal T_*^{-\mathsf T}W\|^{-1}\ne0.
\]
\end{proof}

For the unfolding conditions, let $D_\lambda$ denote the row
vector of partial derivatives with respect to the $m$ parameters,
and set
\[
    c_*=c(0,\lambda_*),\qquad
    c_s^*=\partial_sc(0,\lambda_*).
\]

Let $k\geq2$ be an integer.

\begin{corollary}[Representation-invariant bifurcation conditions]
\label{cor:intrinsic-bifurcation-data}
If $\partial_s^j\psi(0,\lambda_*)=0$ for $0\leq j<k$, then
\begin{equation}
    \partial_s^k\widetilde\psi(0,\lambda_*)
       =c_*\,\partial_s^k\psi(0,\lambda_*).
    \label{eq:invariant-leading-derivative}
\end{equation}
Moreover,
\begin{equation}
    \begin{pmatrix}
      D_\lambda\widetilde\psi\\[1mm]
      D_\lambda\partial_s\widetilde\psi
    \end{pmatrix}_{(0,\lambda_*)}
    =
    \begin{pmatrix}
       c_*&0\\ c_s^*&c_*
    \end{pmatrix}
    \begin{pmatrix}
      D_\lambda\psi\\[1mm]
      D_\lambda\partial_s\psi
    \end{pmatrix}_{(0,\lambda_*)}.
    \label{eq:invariant-unfolding-rows}
\end{equation}
In particular, the rank of the unfolding matrix is independent of
the regular implicit representation.
\end{corollary}

\begin{proof}
Differentiate \eqref{eq:scalar-representation} and use
\eqref{eq:LS-base-derivatives}. The higher state-derivative identity
follows from the product rule when the lower derivatives vanish.
The parameter-derivative identity follows by differentiating once
with respect to $\lambda$ and then with respect to $s$. Its left
factor is invertible because $c_*\ne0$.
\end{proof}

\subsection{Arbitrary-order reduced jet}
\label{subsec:arbitrary-jet}

The derivatives through order $n$ at a point form the $n$-jet.
We compute the coefficients of $\psi(s,\lambda_*)$ recursively,
using the same linear operator at every order. Define
\[
    Y_*(s)=Y(s,\lambda_*),\qquad
    \Xi_n=\frac{d^n}{ds^n}Y_*(s)\bigg|_{s=0}
    \quad(n\geq1).
\]
Equation \eqref{eq:LS-base-derivatives} gives
\begin{equation}
    \Xi_1=V,\qquad
    \Xi_n=\partial_s^nZ(0,\lambda_*)\quad(n\geq2).
    \label{eq:Xi-def}
\end{equation}
Let $\mathfrak P_n$ be the set of partitions of
$\{1,\ldots,n\}$. A partition
$\pi=\{\omega_1,\ldots,\omega_k\}$ consists of $k$ nonempty,
disjoint blocks whose union is $\{1,\ldots,n\}$, and
$|\omega_i|$ is the number of elements in its $i$th block.
For $k\geq2$, $D_X^k\Phi[h_1,\ldots,h_k]$ denotes the symmetric
$k$-linear derivative of $\Phi$ with respect to $X$, evaluated at
$(X_*,\lambda_*)$. Its symmetry makes the ordering of the blocks
immaterial.

For each integer $n\geq2$, define
\begin{equation}
    G_n=
    \sum_{\substack{\pi=\{\omega_1,\ldots,\omega_k\}
                       \in\mathfrak P_n\\ k\geq2}}
    D_X^k\Phi\bigl[\Xi_{|\omega_1|},\ldots,
                        \Xi_{|\omega_k|}\bigr].
    \label{eq:Gn}
\end{equation}
Every term of $G_n$ depends only on
$\Xi_1,\ldots,\Xi_{n-1}$.

\begin{theorem}[Recursive reduced jet]
\label{thm:recursive-jet}
For every $n\geq2$,
\begin{align}
    \Xi_n&=-\mathcal R G_n,
    \label{eq:Xi-recursion}\\
    \partial_s^n\psi(0,\lambda_*)&=\langle W,G_n\rangle.
    \label{eq:psi-n-recursion}
\end{align}
\end{theorem}

\begin{proof}
The higher chain rule, indexed by partitions of
$\{1,\ldots,n\}$, gives
\[
    \frac{d^n}{ds^n}\Phi(Y_*(s),\lambda_*)\bigg|_{s=0}
       =L\Xi_n+G_n.
\]
The term $L\Xi_n$ corresponds to the partition with one block.
Differentiating the range equation $n$ times gives
$L\Xi_n+QG_n=0$, because $QL=L$. For $n\geq2$,
$\Xi_n\in V^\perp$; applying the reduced inverse gives
\eqref{eq:Xi-recursion}. Taking the inner product of the full
derivative with $W$ gives \eqref{eq:psi-n-recursion}, since
$W\transp L=0$.
\end{proof}

\section{Reduced Coefficients and the \texorpdfstring{$+1$}{+1} Bifurcations}
\label{sec:plus-one}

The reduction of Section~\ref{sec:LS} describes the periodic orbits
near a simple multiplier $+1$ by the scalar equation
$\psi(s,\lambda)=0$. Its quadratic and cubic state derivatives
and the leading parameter derivatives determine the fold,
transcritical, pitchfork and cusp bifurcations. We first express
these derivatives in terms of the original implicit equations.
Throughout the next two subsections, derivatives of $\Phi$ are
evaluated at $(X_*,\lambda_*)$ and derivatives of $\psi$ at
$(0,\lambda_*)$. We write
\begin{equation}
    \Phi_i(X,\lambda)=\partial_{\lambda_i}\Phi(X,\lambda)
    \label{eq:param-derivatives}
\end{equation}
for the derivative of the cyclic orbit map with respect to the
$i$th parameter.

\subsection{Quadratic and cubic coefficients}

The recursion in Theorem~\ref{thm:recursive-jet} begins with
\begin{equation}
    Z_{ss}(0,\lambda_*)=\Xi_2
       =-\mathcal R D_X^2\Phi[V,V].
    \label{eq:Zss}
\end{equation}
This correction accounts for the displacement of the noncritical
orbit variables as the critical coordinate $s$ varies.

\begin{proposition}[Quadratic and cubic reduced coefficients]
\label{prop:a2a3}
At a simple multiplier $+1$,
\begin{equation}
    \psi_{ss}(0,\lambda_*)
       =\langle W,D_X^2\Phi[V,V]\rangle,
    \label{eq:psi-ss}
\end{equation}
and
\begin{equation}
    \psi_{sss}(0,\lambda_*)
       =\Bigl\langle W,
          D_X^3\Phi[V,V,V]
          -3D_X^2\Phi\bigl[V,\mathcal R D_X^2\Phi[V,V]\bigr]
         \Bigr\rangle.
    \label{eq:psi-sss}
\end{equation}
\end{proposition}

\begin{proof}
Theorem~\ref{thm:recursive-jet} gives
\[
    G_2=D_X^2\Phi[V,V],\qquad
    G_3=D_X^3\Phi[V,V,V]+3D_X^2\Phi[V,\Xi_2].
\]
Pairing these expressions with $W$ and substituting
\eqref{eq:Zss} proves both formulas.
\end{proof}

The formulas can be expressed directly through the individual
implicit relations. Write $V=(v_0,\ldots,v_{q-1})$,
$W=(w_0,\ldots,w_{q-1})$, and
$Z_{ss}(0,\lambda_*)=(z_0,\ldots,z_{q-1})$. With cyclic indices, set
$h_j=(v_j,v_{j+1})$ and $\widehat h_j=(z_j,z_{j+1})$.
Write $D^kF_j$ for the $k$th derivative in the two state arguments,
with the parameter fixed, evaluated at $(x_j^*,x_{j+1}^*,\lambda_*)$. Then
\[
\begin{aligned}
\psi_{ss}(0,\lambda_*)
 &=\sum_{j=0}^{q-1}\bigl\langle w_j,D^2F_j[h_j,h_j]\bigr\rangle,\\
\psi_{sss}(0,\lambda_*)
 &=\sum_{j=0}^{q-1}\bigl\langle w_j,
     D^3F_j[h_j,h_j,h_j]+3D^2F_j[h_j,\widehat h_j]\bigr\rangle.
\end{aligned}
\]
The blocks $z_j$ are obtained from \eqref{eq:Zss}; only derivatives
of consecutive-state relations and a linear solve are needed.

\subsection{Parameter derivatives}

Differentiating the range equation \eqref{eq:range-equation}
with respect to $\lambda_i$ gives
\begin{equation}
    Z_{\lambda_i}(0,\lambda_*)=-\mathcal R\Phi_i.
    \label{eq:Zlambda}
\end{equation}

\begin{proposition}[Leading unfolding derivatives]
\label{prop:param-derivatives}
For each parameter $\lambda_i$,
\begin{equation}
    \psi_{\lambda_i}(0,\lambda_*)
       =\langle W,\Phi_i\rangle,
    \label{eq:psi-lambda}
\end{equation}
and
\begin{equation}
    \psi_{s\lambda_i}(0,\lambda_*)
       =\Bigl\langle W,
          D_X\Phi_i[V]
          -D_X^2\Phi\bigl[V,\mathcal R\Phi_i\bigr]
         \Bigr\rangle.
    \label{eq:psi-slambda}
\end{equation}
\end{proposition}

\begin{proof}
Differentiate $\psi=\langle W,\Phi(Y,\lambda)\rangle$.
At the base point, $Y_s=V$ and $Y_{\lambda_i}=Z_{\lambda_i}$.
Terms containing $L Y_{\lambda_i}$ or $L Y_{s\lambda_i}$
vanish when paired with $W$. Hence
\[
    \psi_{\lambda_i}=\langle W,\Phi_i\rangle,
    \qquad
    \psi_{s\lambda_i}
       =\langle W,D_X\Phi_i[V]
                  +D_X^2\Phi[V,Z_{\lambda_i}]\rangle.
\]
Substitution of \eqref{eq:Zlambda} completes the proof.
\end{proof}

Two reduced equations are \emph{locally equivalent} if one is obtained
from the other by a smooth change of scalar state coordinate depending
on the parameters, a smooth local change of parameter coordinates,
and multiplication by a smooth nonvanishing function. Both coordinate
changes are local diffeomorphisms and preserve the base point.
Write $\xi$ for the new state coordinate and $\eta$, or
$(\eta_1,\eta_2)$, for the new parameters. The normal forms below
classify bifurcation equations under this equivalence. Their signs
$\pm$ do not distinguish equivalence classes; existence conditions
are also given in the original parameters.

Write $\Sone=\{\zeta\in\C:|\zeta|=1\}$ for the unit circle.
We impose the following spectral isolation condition.

\begin{assumption}[Spectral isolation at $+1$]
\label{ass:isolation-plus}
The multiplier $+1$ is algebraically simple and is the only
multiplier on the unit circle:
\begin{equation}
    \spec(M)\cap\Sone=\{1\}.
    \label{eq:isolation-plus}
\end{equation}
\end{assumption}

\subsection{Fold}

Consider a one-parameter family with parameter $\mu$, centred
at the critical value $\mu=0$.

\begin{theorem}[Fold of a periodic orbit]
\label{thm:fold}
Under Assumption~\ref{ass:isolation-plus}, suppose that
\begin{equation}
    \psi_\mu(0,0)\ne0,
    \qquad
    \psi_{ss}(0,0)\ne0.
    \label{eq:fold-conditions}
\end{equation}
Then the reduced equation is locally equivalent to
\begin{equation}
    \eta\pm\xi^2=0.
    \label{eq:fold-normal}
\end{equation}
\end{theorem}

\begin{proof}
Since $\psi_\mu(0,0)\ne0$, by the implicit function theorem
we can write the solution set as $\mu=b(s)$, where $b(0)=0$.
Differentiating $\psi(s,b(s))=0$ at the origin gives
\[
    b'(0)=0,
    \qquad
    b''(0)=-\frac{\psi_{ss}(0,0)}{\psi_\mu(0,0)}\ne0.
\]
Thus $b(s)=s^2u(s)$ for a smooth function with $u(0)\ne0$.
The change of coordinate $\xi=s\sqrt{|u(s)|}$ brings the
solution curve to $\mu=\pm\xi^2$. The defining equation differs
from $\mu-b(s)$ by a smooth nonvanishing factor, which gives
\eqref{eq:fold-normal}.
\end{proof}

In the original parameter, the fold locus satisfies
\[
    \mu=-\frac{\psi_{ss}(0,0)}{2\psi_\mu(0,0)}s^2+O(s^3).
\]
For sufficiently small nonzero $\mu$, there are two nearby
solutions precisely when
\[
    \mu\,\psi_\mu(0,0)\psi_{ss}(0,0)<0.
\]

\subsection{Persistent branches: transcritical and pitchfork}

For a transcritical or pitchfork bifurcation, a branch of
$q$-periodic orbits continues through the critical parameter.
Let $X_0(\mu)$ be such a smooth branch, with
\[
    \Phi(X_0(\mu),\mu)=0,
    \qquad X_0(0)=X_*.
\]
By Theorem~\ref{thm:LS}, its reduced coordinate is
\[
    s_0(\mu)=\langle V,X_0(\mu)-X_*\rangle,
    \qquad X_0(\mu)=Y(s_0(\mu),\mu).
\]
To place this branch at the origin of the state coordinate, define
\begin{equation}
    g(t,\mu)=\psi(s_0(\mu)+t,\mu).
    \label{eq:persistent-branch}
\end{equation}
Then $g(0,\mu)=0$ for all nearby $\mu$ and $s_0(0)=0$.
At the critical point,
\begin{equation}
\begin{aligned}
    g_{tt}(0,0)&=\psi_{ss}(0,0),\\
    g_{ttt}(0,0)&=\psi_{sss}(0,0),\\
    g_{t\mu}(0,0)&=\psi_{s\mu}(0,0)
                        +s_0'(0)\psi_{ss}(0,0).
\end{aligned}
    \label{eq:centred-derivatives}
\end{equation}
In particular, the mixed derivative governing the bifurcation
is that of the equation centred on the persistent branch.

The identity $g(0,\mu)=0$ gives a smooth factorisation
\begin{equation}
    g(t,\mu)=t\,\chi(t,\mu),
    \qquad
    \chi(t,\mu)=\int_0^1 g_t(\theta t,\mu)\,d\theta.
    \label{eq:persistent-factor}
\end{equation}
At the critical point, $\chi(0,0)=0$ and
\begin{equation}
    \chi_\mu(0,0)=g_{t\mu}(0,0),\qquad
    \chi_t(0,0)=\tfrac12g_{tt}(0,0),\qquad
    \chi_{tt}(0,0)=\tfrac13g_{ttt}(0,0).
    \label{eq:chi-derivatives}
\end{equation}

\begin{theorem}[Transcritical bifurcation]
\label{thm:transcritical}
Under Assumption~\ref{ass:isolation-plus}, suppose that
\begin{equation}
    g_{t\mu}(0,0)\ne0,
    \qquad
    g_{tt}(0,0)\ne0.
    \label{eq:transcritical-conditions}
\end{equation}
Then the reduced equation is locally equivalent to
\begin{equation}
    \xi(\eta\pm\xi)=0.
    \label{eq:transcritical-normal}
\end{equation}
\end{theorem}

\begin{proof}
By \eqref{eq:chi-derivatives}, $\chi_\mu(0,0)$ and
$\chi_t(0,0)$ are nonzero. Hence $\chi(t,\mu)=0$ defines a
smooth curve $\mu=b(t)$ with
\[
    b(0)=0,\qquad
    b'(0)=-\frac{\chi_t(0,0)}{\chi_\mu(0,0)}\ne0.
\]
Writing $\chi(t,\mu)=u(t,\mu)(\mu-b(t))$ with $u(0,0)\ne0$,
we obtain $g=t\,u(t,\mu)(\mu-b(t))$.
The coordinate $\xi=b(t)$ is invertible near the origin; after
absorbing the nonzero factors, the equation becomes
\eqref{eq:transcritical-normal}.
\end{proof}

\begin{theorem}[Pitchfork bifurcation]
\label{thm:pitchfork}
Under Assumption~\ref{ass:isolation-plus}, suppose that
\begin{equation}
    g_{t\mu}(0,0)\ne0,
    \qquad
    g_{tt}(0,0)=0,
    \qquad
    g_{ttt}(0,0)\ne0.
    \label{eq:pitchfork-conditions}
\end{equation}
Then the reduced equation is locally equivalent to
\begin{equation}
    \xi(\eta\pm\xi^2)=0.
    \label{eq:pitchfork-normal}
\end{equation}
\end{theorem}

\begin{proof}
By \eqref{eq:chi-derivatives},
$\chi_\mu(0,0)\ne0$, $\chi_t(0,0)=0$ and
$\chi_{tt}(0,0)\ne0$. The curve $\chi(t,\mu)=0$ can therefore
be written $\mu=b(t)$, with
\[
    b(0)=b'(0)=0,
    \qquad
    b''(0)=-\frac{\chi_{tt}(0,0)}{\chi_\mu(0,0)}\ne0.
\]
Write $b(t)=t^2u(t)$, where $u(0)\ne0$, and set
$\xi=t\sqrt{|u(t)|}$. Factoring $\chi$ as in the proof of
Theorem~\ref{thm:transcritical} gives
\eqref{eq:pitchfork-normal}.
\end{proof}

In the original parameter, the two branches away from the
persistent one satisfy
\[
    t^2=-6\frac{g_{t\mu}(0,0)}{g_{ttt}(0,0)}\mu+o(|\mu|).
\]
They exist for sufficiently small nonzero $\mu$ precisely when
$\mu\,g_{t\mu}(0,0)g_{ttt}(0,0)<0$.

\subsection{Cusp}

For the cusp, take two parameters $\lambda=(\lambda_1,\lambda_2)$
and place the critical parameter at $\lambda=0$. The first
nonzero state term is cubic, and the two parameters independently
vary the constant and linear coefficients of the reduced equation.

\begin{theorem}[Cusp of a periodic orbit]
\label{thm:cusp}
Under Assumption~\ref{ass:isolation-plus}, suppose that
\begin{equation}
    \psi_{ss}(0,0)=0,
    \qquad
    \psi_{sss}(0,0)\ne0.
    \label{eq:cusp-degeneracy}
\end{equation}
If
\begin{equation}
    \rank
    \begin{pmatrix}
      \psi_{\lambda_1}(0,0) & \psi_{\lambda_2}(0,0)\\
      \psi_{s\lambda_1}(0,0) & \psi_{s\lambda_2}(0,0)
    \end{pmatrix}=2,
    \label{eq:cusp-rank}
\end{equation}
then the reduced equation is locally equivalent to
\begin{equation}
    \xi^3+\eta_1\xi+\eta_2=0,
    \label{eq:cusp-normal}
\end{equation}
\end{theorem}

\begin{proof}
At the origin, $\psi=\psi_s=\psi_{ss}=0$, while
$\psi_{sss}\ne0$. The smooth preparation theorem
\citep{GolubitskySchaeffer1985} therefore gives
\[
    \psi(s,\lambda)=u(s,\lambda)
       \bigl(s^3+c_2(\lambda)s^2
                 +c_1(\lambda)s+c_0(\lambda)\bigr),
\]
where $u(0,0)\ne0$ and $c_j(0)=0$ for $j=0,1,2$.
The rank condition \eqref{eq:cusp-rank} implies that
$\lambda\mapsto(c_0(\lambda),c_1(\lambda))$ has invertible
derivative at the origin: differentiating the prepared expression
at $s=0$ gives
\[
    \begin{pmatrix}
      D_\lambda\psi\\ D_\lambda\psi_s
    \end{pmatrix}_{(0,0)}
    =\begin{pmatrix}
       u(0,0)&0\\ u_s(0,0)&u(0,0)
      \end{pmatrix}
      \begin{pmatrix}
        D_\lambda c_0(0)\\ D_\lambda c_1(0)
      \end{pmatrix}.
\]
The translation $s=\xi-c_2(\lambda)/3$ removes the quadratic
term. The resulting coefficients are
\[
\begin{aligned}
    \eta_1(\lambda)
       &=c_1(\lambda)-\tfrac13c_2(\lambda)^2,\\
    \eta_2(\lambda)
       &=c_0(\lambda)-\tfrac13c_1(\lambda)c_2(\lambda)
                         +\tfrac{2}{27}c_2(\lambda)^3.
\end{aligned}
\]
Their parameter derivatives at the origin are those of $c_1$
and $c_0$, respectively. The rank condition makes
$(\eta_1,\eta_2)$ valid local parameters; division by $u$
then gives \eqref{eq:cusp-normal}.
\end{proof}

\section{Period Doubling and the Doubled Cyclic System}
\label{sec:doubling}

A perturbation associated with a multiplier $-1$ changes sign after
one period. It closes after two periods, so period doubling is naturally
studied by writing two consecutive copies of the periodic-orbit
equations. Interchanging the copies leaves those equations unchanged.
This exchange symmetry is the basis of the Lyapunov--Schmidt approach
to degenerate period doubling developed for explicit maps by Peckham
and Kevrekidis~\citeyearpar{PeckhamKevrekidis1991}. The exchange symmetry is a particular instance of the
equivariant bifurcation setting discussed in
\citet{GolubitskyStewartSchaeffer1988}.

Fix a periodic solution $(X_*,0)$.

\begin{assumption}[Spectral isolation at $-1$]
\label{ass:isolation-minus}
At $\lambda=0$,
\begin{equation}
    \spec(M)\cap\Sone=\{-1\},
    \qquad -1\text{ is algebraically simple}.
    \label{eq:isolation-minus}
\end{equation}
\end{assumption}

Since $+1$ is excluded, $D_X\Phi(X_*,0)$ is invertible.
The implicit function theorem gives a unique smooth branch
$X_0(\lambda)$ with $\Phi(X_0(\lambda),\lambda)=0$ and $X_0(0)=X_*$.

For the doubled orbit space, write
\[
    \widehat X=(X^{(0)},X^{(1)}),\qquad
    X^{(\varepsilon)}=(x_0^{(\varepsilon)},\ldots,x_{q-1}^{(\varepsilon)})
       \in\prod_{j=0}^{q-1}U_j,\quad \varepsilon=0,1.
\]
The doubled cyclic orbit map
\[
    \widehat\Phi:
    \left(\prod_{j=0}^{q-1}U_j\right)^2\times\Lambda
    \longrightarrow\R^{2dq}
\]
has components
\[
\begin{aligned}
 (\widehat\Phi(\widehat X,\lambda))^{(\varepsilon)}_j
     &=F_j(x_j^{(\varepsilon)},x_{j+1}^{(\varepsilon)},\lambda),
           &&0\leq j<q-1,\\
 (\widehat\Phi(\widehat X,\lambda))^{(\varepsilon)}_{q-1}
     &=F_{q-1}(x_{q-1}^{(\varepsilon)},x_0^{(1-\varepsilon)},\lambda),
           &&\varepsilon=0,1.
\end{aligned}
\]
The second line joins each block to the next. In particular,
\[
    \widehat X_0(\lambda)
       =(X_0(\lambda),X_0(\lambda))
\]
is a solution. A solution with unequal blocks represents a two-cycle
of the $q$-step evolution. If the original orbit has least period
$q$, nearby unequal-block solutions have least period $2q$;
otherwise the conclusion concerns the $q$-step return.
Define the exchange operator
\begin{equation}
    \mathcal S(X^{(0)},X^{(1)})=(X^{(1)},X^{(0)}).
    \label{eq:block-exchange}
\end{equation}
It is orthogonal, $\mathcal S^2=I$, and the periodicity of the
defining relations gives
\begin{equation}
    \widehat\Phi(\mathcal S\widehat X,\lambda)
      =\mathcal S\widehat\Phi(\widehat X,\lambda).
    \label{eq:doubled-equivariance}
\end{equation}

At $\widehat X_*=(X_*,X_*)$, let
$\widehat L=D_{\widehat X}\widehat\Phi(\widehat X_*,0)$.
The exchange operator splits the doubled space into its symmetric
and antisymmetric subspaces,
\[
    \mathcal E^+=\{(U,U):U\in\R^{dq}\},\qquad
    \mathcal E^-=\{(U,-U):U\in\R^{dq}\}.
\]
For $U=(u_0,\ldots,u_{q-1})\in\R^{dq}$, set
$\iota_\pm U=(U,\pm U)$.

\begin{proposition}[Symmetric--antisymmetric splitting]
\label{prop:doubled-splitting}
Both $\mathcal E^+$ and $\mathcal E^-$ are invariant under
$\widehat L$, and
\begin{equation}
    \widehat L\,\iota_+U=\iota_+(L_1U),\qquad
    \widehat L\,\iota_-U=\iota_-(L_{-1}U).
    \label{eq:doubled-restrictions}
\end{equation}
Consequently,
\begin{equation}
    \Ker\widehat L
      =\iota_+(\Ker L_1)\oplus\iota_-(\Ker L_{-1}).
    \label{eq:doubled-kernel}
\end{equation}
Under Assumption~\ref{ass:isolation-minus}, $\widehat L$ has a
one-dimensional kernel contained in $\mathcal E^-$.
\end{proposition}

\begin{proof}
For the first $q-1$ equations in each block, the linearisation
acts as $A_j u_j+B_j u_{j+1}$. In the last equation the next state
belongs to the other block. Substitution of $(U,U)$ gives the closing
term $B_{q-1}u_0$, whereas $(U,-U)$ gives $-B_{q-1}u_0$ in the first
block and its negative in the second. These are exactly the closing
terms of $L_1$ and $L_{-1}$, proving
\eqref{eq:doubled-restrictions}. Since
$\R^{2dq}=\mathcal E^+\oplus\mathcal E^-$, the kernel splits as
stated. By Theorem~\ref{thm:twisted-monodromy}, the spectral
assumption makes $L_1$ invertible and gives
$\dim\Ker L_{-1}=1$.
\end{proof}

\subsection{Lyapunov--Schmidt reduction of the doubled system}

Choose unit vectors $V_-\in\Ker L_{-1}$ and
$W_-\in\Ker L_{-1}\transp$, and define
\[
    \widehat V=\frac{1}{\sqrt2}\iota_-V_-,
    \qquad
    \widehat W=\frac{1}{\sqrt2}\iota_-W_-.
\]
Proposition~\ref{prop:doubled-splitting} shows that these vectors
span the right and left kernels of $\widehat L$. In particular,
$\mathcal S\widehat V=-\widehat V$ and
$\mathcal S\widehat W=-\widehat W$. Let
$\widehat Q=I-\widehat W\widehat W\transp$.

We apply the reduction of Section~\ref{sec:LS} with the critical
coordinate centred on the repeated branch. For each $(s,\lambda)$
near $(0,0)$, the complementary equation
\[
    \widehat Q\,\widehat\Phi\bigl(
      \widehat X_0(\lambda)+s\widehat V+\widehat Z,\lambda
      \bigr)=0,
    \qquad \widehat Z\in\widehat V^\perp,
\]
has a unique smooth solution $\widehat Z=\widehat Z(s,\lambda)$.
Since $\widehat X_0(\lambda)$ is a solution,
$\widehat Z(0,\lambda)=0$. The remaining scalar equation is
\[
    \widehat\psi(s,\lambda)
      =\left\langle\widehat W,
        \widehat\Phi\bigl(\widehat X_0(\lambda)
          +s\widehat V+\widehat Z(s,\lambda),\lambda\bigr)
        \right\rangle=0.
\]
Its solutions correspond exactly to the nearby solutions of the
doubled cyclic system.

Define the smooth coefficient functions
\begin{equation}
    \beta_1(\lambda)=\partial_s\widehat\psi(0,\lambda),
    \quad
    \beta_3(\lambda)=\frac{1}{3!}\partial_s^3\widehat\psi(0,\lambda),
    \quad
    \beta_5(\lambda)=\frac{1}{5!}\partial_s^5\widehat\psi(0,\lambda).
    \label{eq:beta35}
\end{equation}

\begin{theorem}[Odd reduction]
\label{thm:odd-reduction}
The reduced equation satisfies
\begin{equation}
    \widehat\psi(-s,\lambda)=-\widehat\psi(s,\lambda).
    \label{eq:oddness}
\end{equation}
\end{theorem}

Consequently,
\begin{equation}
    \widehat\psi(s,\lambda)
     =s\left(\beta_1(\lambda)+\beta_3(\lambda)s^2
                 +\beta_5(\lambda)s^4+O(s^6)\right),
    \qquad \beta_1(0)=0,
    \label{eq:odd-reduced}
\end{equation}

At $\lambda=0$, the coefficients follow from the recursion in
Theorem~\ref{thm:recursive-jet} applied to $\widehat\Phi$.

\begin{proof}
The repeated branch is fixed by $\mathcal S$. Applying
\eqref{eq:doubled-equivariance} to the complementary equation and
using $\mathcal S\widehat V=-\widehat V$ shows that
$\mathcal S\widehat Z(s,\lambda)$ solves the complementary
equation at $-s$. Its uniqueness gives
\[
    \widehat Z(-s,\lambda)=\mathcal S\widehat Z(s,\lambda).
\]
The orthogonality of $\mathcal S$ and
$\mathcal S\widehat W=-\widehat W$ now give
$\widehat\psi(-s,\lambda)=-\widehat\psi(s,\lambda)$.
The expansion follows from smoothness and oddness. Finally,
$\beta_1(0)=0$ because $\widehat V\in\Ker\widehat L$.
\end{proof}

Thus the solution $s=0$ is the repeated $q$-periodic orbit; every
nonzero root occurs together with its negative. These two roots
represent the same two-cycle with the two $q$-blocks interchanged.

The linear coefficient also measures the crossing of the multiplier
through $-1$. Let $M(\lambda)$ be the monodromy matrix along the
branch $X_0(\lambda)$, and let $\kappa(\lambda)$ be its simple real
multiplier with $\kappa(0)=-1$. Write $L_\zeta(\lambda)$ for the
twisted operator along this branch, and put
\[
    \gamma
       =\left\langle W_-,
          \left.\partial_\zeta L_\zeta(0)\right|_{\zeta=-1}V_-
        \right\rangle.
\]

\begin{proposition}[Multiplier crossing]
\label{prop:flip-crossing}
Under Assumption~\ref{ass:isolation-minus}, $\gamma\ne0$ and
\begin{equation}
    D_\lambda\beta_1(0)=-\gamma\,D_\lambda\kappa(0).
    \label{eq:flip-crossing}
\end{equation}
In particular, $\partial_\mu\beta_1(0)\ne0$ in a one-parameter
family if and only if $\kappa'(0)\ne0$.
\end{proposition}

\begin{proof}
Differentiate the reduced equation with respect to $s$ at
$s=0$. At the critical point, $\widehat Z_s(0,0)=0$, and terms
containing $\widehat L$ vanish when paired with $\widehat W$.
Differentiation along the branch therefore gives, for each
parameter $\lambda_i$,
\[
    \partial_{\lambda_i}\beta_1(0)
      =\left\langle W_-,
        \partial_{\lambda_i}L_{-1}(0)\,V_-\right\rangle,
\]
where the parameter derivative of $L_{-1}$ includes the variation
of $X_0(\lambda)$. On the other hand,
$L_{\kappa(\lambda)}(\lambda)$ has a nonzero kernel.
Differentiating its kernel equation at $\lambda=0$ and pairing
with $W_-$ gives
\[
    \left\langle W_-,
       \partial_{\lambda_i}L_{-1}(0)V_-\right\rangle
       +\gamma\,\partial_{\lambda_i}\kappa(0)=0.
\]
To prove $\gamma\ne0$, the determinant identity
\eqref{eq:det-identity} and the algebraic simplicity of $-1$ give
\[
    \left.\partial_\zeta\det L_\zeta(0)\right|_{\zeta=-1}\ne0.
\]
Since $L_{-1}(0)$ has rank $dq-1$, its adjugate is
$aV_-W_-\transp$ for some $a\ne0$. The derivative formula for a
determinant therefore gives
\[
    \left.\partial_\zeta\det L_\zeta(0)\right|_{\zeta=-1}
       =a\left\langle W_-,
             \left.\partial_\zeta L_\zeta(0)\right|_{\zeta=-1}
              V_-\right\rangle
       =a\gamma.
\]
Thus $\gamma\ne0$.
\end{proof}

\subsection{Flip and degenerate flip}

Since $\widehat\psi$ is smooth and odd in $s$, Hadamard's lemma
gives $\widehat\psi(s,\lambda)=s\,h(s,\lambda)$, where $h$ is smooth
and even in $s$. The smooth even-function lemma, valid also with
parameters, gives a smooth function $\varphi(t,\lambda)$ near $(0,0)$ such that
$h(s,\lambda)=\varphi(s^2,\lambda)$. Indeed, for $t>0$ set
$\varphi(t,\lambda)=h(\sqrt{t},\lambda)$. The operator
$(2s)^{-1}\partial_s$ maps smooth even functions to smooth even functions locally,
by Hadamard's lemma. Iteration shows that all derivatives of
$\varphi$ in $t$ and $\lambda$ extend continuously to $t=0$;
a smooth extension across $t=0$ completes the construction.
Consequently, $\widehat\psi(s,\lambda)=s\,\varphi(s^2,\lambda)$. Thus
\[
    \varphi(0,\lambda)=\beta_1(\lambda),\qquad
    \partial_t\varphi(0,\lambda)=\beta_3(\lambda),\qquad
    \tfrac12\partial_t^2\varphi(0,\lambda)=\beta_5(\lambda),
\]
where $t$ denotes the argument of $\varphi$. The root $s=0$ represents
the persistent orbit; nonzero roots are determined by
$\varphi(t,\lambda)=0$ with $t=s^2>0$.

For the ordinary flip, use one parameter $\mu$, centred at zero.

\begin{theorem}[Flip]
\label{thm:flip}
Under Assumption~\ref{ass:isolation-minus}, if
\begin{equation}
    \partial_\mu\beta_1(0)\ne0,\qquad \beta_3(0)\ne0,
    \label{eq:flip-conditions}
\end{equation}
then the doubled reduced equation is locally equivalent to
\begin{equation}
    \xi(\eta\pm\xi^2)=0.
    \label{eq:flip-normal}
\end{equation}
\end{theorem}

\begin{proof}
For the ordinary flip, $\partial_t\varphi(0,0)=\beta_3(0)\ne0$.
The implicit function theorem gives a smooth function
$t=t_*(\mu)$ describing the zeros of $\varphi$, with
\[
    t_*(0)=0,\qquad
    t_*'(0)=-\frac{\partial_\mu\beta_1(0)}{\beta_3(0)}\ne0.
\]
By the fundamental theorem of calculus,
$\varphi(t,\mu)=u(t,\mu)(t-t_*(\mu))$, where $u(0,0)\ne0$.
The change of parameter $\eta=-t_*(\mu)$ and division by
the nonvanishing factor $u$ give \eqref{eq:flip-normal}.
\end{proof}

For the degenerate flip, use two parameters
$\lambda=(\lambda_1,\lambda_2)$, centred at zero.

\begin{theorem}[Degenerate flip]
\label{thm:flip-gpd}
Under Assumption~\ref{ass:isolation-minus}, if
\begin{equation}
    \beta_3(0)=0,\qquad \beta_5(0)\ne0,
    \label{eq:gpd-degeneracy}
\end{equation}
and
\begin{equation}
    \rank D_\lambda(\beta_1,\beta_3)(0)=2,
    \label{eq:gpd-rank}
\end{equation}
then the doubled reduced equation is locally equivalent to
\begin{equation}
    \xi(\eta_1+\eta_2\xi^2\pm\xi^4)=0.
    \label{eq:gpd-normal}
\end{equation}
\end{theorem}

\begin{proof}
For the degenerate flip, $\varphi(0,0)=\partial_t\varphi(0,0)=0$
and $\partial_t^2\varphi(0,0)=2\beta_5(0)\ne0$.
Smooth preparation~\citep{GolubitskySchaeffer1985} gives
\[
    \varphi(t,\lambda)
      =u(t,\lambda)\bigl(t^2+a_1(\lambda)t+a_0(\lambda)\bigr),
    \qquad u(0,0)\ne0,
\]
with $a_0(0)=a_1(0)=0$. Differentiating at $t=0$, one finds
\[
    \begin{pmatrix}D_\lambda\beta_1\\D_\lambda\beta_3\end{pmatrix}_{\lambda=0}
      =
    \begin{pmatrix}u_*&0\\u_{t,*}&u_*\end{pmatrix}
    \begin{pmatrix}D_\lambda a_0\\D_\lambda a_1\end{pmatrix}_{\lambda=0},
    \qquad u_*=u(0,0),\quad u_{t,*}=\partial_t u(0,0).
\]
The triangular matrix is invertible, so
\eqref{eq:gpd-rank} makes $(a_0,a_1)$ valid local parameter
coordinates. Taking $\eta_1=a_0$, $\eta_2=a_1$ and dividing by
$u$ gives \eqref{eq:gpd-normal}. In this form, the nonzero
solutions correspond to the positive roots of a quadratic
polynomial in $s^2$.
\end{proof}

For the ordinary flip, the nonzero roots of the doubled equation
satisfy, in the original parameter,
\[
    s^2=-\frac{\partial_\mu\beta_1(0)}{\beta_3(0)}\mu
          +O(\mu^2).
\]
Consequently, a nearby period-doubled orbit exists for sufficiently
small nonzero $\mu$ precisely when
\[
    \mu\,\partial_\mu\beta_1(0)\beta_3(0)<0,
    \quad\text{equivalently}\quad
    \mu\,\gamma\,\kappa'(0)\beta_3(0)>0
\]
by Proposition~\ref{prop:flip-crossing}.

Applying the change of representation to each equation of the doubled system gives
$\widetilde{\widehat\Phi}=\widehat{\mathcal T}\widehat\Phi$.
We denote its reduction and coefficients by $\widetilde{\widehat\psi}$ and $\widetilde\beta_j$.
Centre both reductions on the same repeated branch, retain $\widehat V$,
and choose $\widetilde{\widehat W}$ by \eqref{eq:transformed-left-vector}
with $\mathcal T_*$ and $W$ replaced by their doubled counterparts.
Let $c_{\mathrm d}$ be the comparison factor constructed above
Theorem~\ref{thm:reduced-representation}, applied to these doubled
reductions, and set
\[
    c(s,\lambda)=\tfrac12\bigl(c_{\mathrm d}(s,\lambda)
                                +c_{\mathrm d}(-s,\lambda)\bigr).
\]
A smooth nonvanishing function $c(s,\lambda)$ is an \emph{even unit}
if $c(-s,\lambda)=c(s,\lambda)$. For such a function, write
$c_0(\lambda)=c(0,\lambda)$, $c_2(\lambda)=\tfrac12\partial_s^2c(0,\lambda)$,
and $c_j^*=c_j(0)$ for $j=0,2$.
\begin{corollary}[Representation-invariant flip conditions]
\label{cor:flip-representation}
Under Assumption~\ref{ass:isolation-minus}, $c$ is an even unit,
$\widetilde{\widehat\psi}=c\widehat\psi$ and
$\widetilde\beta_3(0)=c_0^*\beta_3(0)$. If $\beta_3(0)=0$, then
\[
 \widetilde\beta_5(0)=c_0^*\beta_5(0),\qquad
 D_\lambda\binom{\widetilde\beta_1}{\widetilde\beta_3}(0)
 =\begin{pmatrix}c_0^*&0\\c_2^*&c_0^*\end{pmatrix}
 D_\lambda\binom{\beta_1}{\beta_3}(0).
\]
Thus the flip and degenerate-flip nondegeneracy and unfolding conditions
are independent of the forward-regular representation.
\end{corollary}
\begin{proof}
The matrix $\widehat{\mathcal T}_*$ commutes with $\mathcal S$,
so the transported left vector remains antisymmetric.
Apply Theorem~\ref{thm:reduced-representation} to the doubled system,
whose multiplier $+1$ is simple, after translating the common branch to the origin. Both reductions
are odd, so symmetrising $c_{\mathrm d}$ preserves the identity
and nonvanishing at the origin. Comparing coefficients gives
$\widetilde\beta_1=c_0\beta_1$ and $\widetilde\beta_3=c_0\beta_3+c_2\beta_1$.
The conclusions follow by differentiation and comparison of the fifth-order terms,
using $\beta_1(0)=0$ and, for the degenerate flip, $\beta_3(0)=0$.
\end{proof}

\subsection{The scalar Schwarzian criterion}

For a scalar relation, the locally defined $q$-step return map
$f=f_{q-1}\circ\cdots\circ f_0$ from Section~\ref{sec:cyclic}
has a fixed point at $x_0^*$ and $f'(x_0^*)=-1$. The Schwarzian
derivative of a $C^3$ map at a point where its derivative is
nonzero is
\[
    \operatorname{Sch}f(x)
       =\frac{f'''(x)}{f'(x)}
          -\frac32\left(\frac{f''(x)}{f'(x)}\right)^2.
\]
The connection between the Schwarzian and scalar period doubling
is classical~\citep{Allwright1978,Singer1978}.

Write $\Id(x)=x$. In the local coordinate $y=x-x_0^*$, define the Taylor coefficients
$r_2,\ldots,r_5$ by
\[
    f(x_0^*+y)-x_0^*=-y+r_2y^2+r_3y^3+r_4y^4+r_5y^5+O(y^6).
\]

\begin{corollary}[Scalar flip coefficients]
\label{cor:flip-schwarzian}
Suppose $d=1$ and Assumption~\ref{ass:isolation-minus} holds.
Then the cubic coefficient in the doubled reduction is nonzero
exactly when $\operatorname{Sch}f(x_0^*)\ne0$. If $\operatorname{Sch}f(x_0^*)=0$, then $r_3=-r_2^2$ and
\begin{equation}
    \beta_5(0)\ne0
       \quad\Longleftrightarrow\quad
    r_5+3r_2r_4-2r_2^4\ne0.
    \label{eq:schwarzian-quintic}
\end{equation}
\end{corollary}

\begin{proof}
Direct substitution gives
\[
    \operatorname{Sch}f(x_0^*)=-6(r_2^2+r_3),
    \qquad
    f^2(x_0^*+y)-x_0^*-y
       =-2(r_2^2+r_3)y^3+O(y^4).
\]
Thus the cubic coefficient of $f^2-\Id$ is
$\operatorname{Sch}f(x_0^*)/3$. When $r_3=-r_2^2$, the fourth-order term
also vanishes, and expansion to fifth order gives
\[
    f^2(x_0^*+y)-x_0^*-y
       =-2(r_5+3r_2r_4-2r_2^4)y^5+O(y^6).
\]
Apply Proposition~\ref{prop:return-equation} to the doubled cyclic
system. Its return map is $f^2$, whose multiplier $+1$ is simple.
Thus $\widehat\psi(\cdot,0)$ and $f^2-\Id$ differ by a smooth local
change of state coordinate and a smooth nonvanishing factor. These
changes preserve the order of the first nonzero derivative, proving
both criteria.
\end{proof}

\section{Cyclic Invariance and Periodic Non-Autonomous Systems}
\label{sec:cyclic-invariance}

A periodic orbit can be indexed from any of its states. In an autonomous
system this only changes the order in which the orbit is written. In a
periodically non-autonomous system the defining relations must be shifted
by the same number of steps. We show that this relabelling leaves the
Lyapunov--Schmidt equations unchanged.
We extend the cyclic invariance established in \citet{Oliveira2016}
for explicit $p$-periodic non-autonomous systems to
regular implicit systems in arbitrary finite dimension, including
period doubling.

For $k\in\{0,\ldots,q-1\}$, let $\rho_k:\R^{dq}\to\R^{dq}$ be the
cyclic permutation of the $q$ state blocks,
\begin{equation}
    (\rho_kX)_j=x_{j+k},
    \qquad j\in\mathbb Z/q\mathbb Z,
    \label{eq:rho-def}
\end{equation}
where indices are read modulo $q$. Shift both the state domains and
the defining relations by setting
\begin{equation}
    U_j^{[k]}=U_{j+k},
    \qquad F_j^{[k]}=F_{j+k}.
    \label{eq:shifted-system}
\end{equation}
Write $\Phi^{[k]}$ for the cyclic orbit map of this shifted system.
Thus $\rho_k$ maps the original orbit domain onto the shifted one.
At the periodic solution $(X_*,\lambda_*)$, set
$L=L_1=D_X\Phi(X_*,\lambda_*)$ and put
\[
    X_*^{[k]}=\rho_kX_*,
    \qquad
    L^{[k]}=D_X\Phi^{[k]}(X_*^{[k]},\lambda_*).
\]
Let $M^{[k]}$ denote the monodromy matrix at the shifted orbit.

\begin{proposition}[Cyclic covariance]
\label{prop:cyclic-covariance}
For every $k$ and every $(X,\lambda)$ in the domain of $\Phi$,
\begin{equation}
    \Phi^{[k]}(\rho_kX,\lambda)
       =\rho_k\Phi(X,\lambda).
    \label{eq:cyclic-covariance}
\end{equation}
In particular,
\begin{equation}
    L^{[k]}=\rho_kL\rho_k^{-1}.
    \label{eq:L-cyclic}
\end{equation}
The monodromy matrices have the same characteristic polynomial:
\begin{equation}
    \det(\zeta I-M^{[k]})=\det(\zeta I-M),
    \qquad \zeta\in\C.
    \label{eq:cyclic-spectrum}
\end{equation}
\end{proposition}

\begin{proof}
The $j$th component of the left-hand side is
$F_{j+k}(x_{j+k},x_{j+k+1},\lambda)$, the $(j+k)$th component of
$\Phi(X,\lambda)$. This proves \eqref{eq:cyclic-covariance}.
Differentiation with respect to $X$ at $X_*$ gives
$L^{[k]}\rho_k=\rho_kL$, and hence \eqref{eq:L-cyclic}.
The forward-step derivatives of the shifted system are
$J_k,J_{k+1},\ldots,J_{k+q-1}$, with subscripts modulo $q$.
Their ordered product is a cyclic rearrangement of
$M=J_{q-1}\cdots J_0$. Successive applications of
$\det(\zeta I-AB)=\det(\zeta I-BA)$ for square matrices
prove \eqref{eq:cyclic-spectrum}.
\end{proof}

\subsection{The multiplier \texorpdfstring{$+1$}{+1}}

Suppose $+1$ is an algebraically simple multiplier. Use the unit kernel
vectors $V,W$ of Section~\ref{sec:LS}. For the shifted problem choose
\begin{equation}
    V^{[k]}=\rho_kV,
    \qquad W^{[k]}=\rho_kW.
    \label{eq:shifted-VW}
\end{equation}
These vectors have unit norm and span the right and left kernels of
$L^{[k]}$. Denote by $\mathcal R^{[k]}$, $Z^{[k]}$ and $\psi^{[k]}$
the reduced inverse, complementary solution and scalar equation
constructed from them as in Section~\ref{sec:LS}.

\begin{theorem}[Cyclic invariance at $+1$]
\label{thm:cyclic-invariance}
For an algebraically simple multiplier $+1$ and the kernel vectors
\eqref{eq:shifted-VW}, the following identities hold for every $k$:
\begin{equation}
    \mathcal R^{[k]}
       =\rho_k\mathcal R\rho_k^{-1},
    \qquad
    Z^{[k]}(s,\lambda)=\rho_kZ(s,\lambda),
    \label{eq:shifted-RZ}
\end{equation}
and
\begin{equation}
    \psi^{[k]}(s,\lambda)=\psi(s,\lambda)
    \label{eq:psi-cyclic-invariant}
\end{equation}
for all $(s,\lambda)$ near $(0,\lambda_*)$.
\end{theorem}

\begin{proof}
The permutation $\rho_k$ is orthogonal. By \eqref{eq:L-cyclic},
it carries $\Ker L$ and $\Ker L\transp$ onto the corresponding
kernels of $L^{[k]}$ and $(L^{[k]})\transp$. It therefore carries
their orthogonal complements onto those of the shifted problem.
In particular,
\[
    Q^{[k]}=I-W^{[k]}(W^{[k]})\transp
            =\rho_kQ\rho_k^{-1}.
\]
Conjugating the restriction of $L$ to $V^\perp$ and using the
definition \eqref{eq:reduced-inverse} proves the first identity
in \eqref{eq:shifted-RZ}.

The vector $\rho_kZ(s,\lambda)$ belongs to $(V^{[k]})^\perp$.
Moreover, \eqref{eq:cyclic-covariance} gives
\[
    Q^{[k]}\Phi^{[k]}
      \bigl(X_*^{[k]}+sV^{[k]}+\rho_kZ(s,\lambda),\lambda\bigr)
       =\rho_kQ\Phi\bigl(Y(s,\lambda),\lambda\bigr)=0.
\]
Uniqueness of the complementary solution proves the second identity
in \eqref{eq:shifted-RZ}. The scalar equations then satisfy
\[
\begin{aligned}
    \psi^{[k]}(s,\lambda)
    &=\bigl\langle\rho_kW,
      \Phi^{[k]}(\rho_kY(s,\lambda),\lambda)\bigr\rangle\\
    &=\bigl\langle\rho_kW,\rho_k\Phi(Y(s,\lambda),\lambda)\bigr\rangle
      =\psi(s,\lambda),
\end{aligned}
\]
by orthogonality of $\rho_k$.
\end{proof}

Thus all state and parameter derivatives of the reduced equation are
unchanged. In particular, the fold and cusp conditions of
Section~\ref{sec:plus-one} are invariant. If a branch of
periodic orbits persists, the coordinate of the shifted branch is
\[
    \bigl\langle V^{[k]},\rho_kX_0(\mu)-X_*^{[k]}\bigr\rangle
     =\bigl\langle V,X_0(\mu)-X_*\bigr\rangle=s_0(\mu).
\]
Hence the centred equation \eqref{eq:persistent-branch} is also
unchanged, giving the same transcritical and pitchfork conditions.

\subsection{The multiplier \texorpdfstring{$-1$}{-1}}

For a flip, we use the doubled system of Section~\ref{sec:doubling}.
Write its $2q$ consecutive state blocks as
$\widehat X=(\widehat x_0,\ldots,\widehat x_{2q-1})$,
where $\widehat x_{\varepsilon q+j}=x_j^{(\varepsilon)}$ for $\varepsilon=0,1$ and $0\leq j<q$.
The corresponding cyclic permutation acts on all $2q$ blocks,
with indices taken modulo $2q$:
\begin{equation}
    (\widehat\rho_k\widehat X)_\ell=\widehat x_{\ell+k},
    \qquad \ell\in\mathbb Z/(2q)\mathbb Z.
    \label{eq:doubled-cyclic-shift}
\end{equation}
This permutation crosses the boundary between the two copies of the
orbit. The exchange operator $\mathcal S$ of
\eqref{eq:block-exchange} is a shift by $q$ blocks. It therefore
commutes with $\widehat\rho_k$, which preserves the antisymmetric
subspace containing the critical direction.

Let $\widehat\Phi^{[k]}$ be the doubled cyclic map formed from the
shifted relations $F_j^{[k]}$. Along the repeated branch, use
\[
    \widehat X_0^{[k]}(\lambda)
        =\widehat\rho_k\widehat X_0(\lambda),
    \qquad
    \widehat V^{[k]}=\widehat\rho_k\widehat V,
    \qquad
    \widehat W^{[k]}=\widehat\rho_k\widehat W
\]
Define
$\widehat L^{[k]}=D_{\widehat X}\widehat\Phi^{[k]}
(\widehat X_0^{[k]}(0),0)$, and write
$\widehat Z^{[k]}$ and $\widehat\psi^{[k]}$ for the
complementary solution and reduced equation of the shifted problem.

\begin{theorem}[Cyclic invariance at $-1$]
\label{thm:cyclic-invariance-flip}
Under Assumption~\ref{ass:isolation-minus}, for every $k$,
\begin{equation}
    \widehat\psi^{[k]}(s,\lambda)
      =\widehat\psi(s,\lambda).
    \label{eq:flip-cyclic-invariant}
\end{equation}
Consequently, the coefficients $\beta_1$, $\beta_3$ and $\beta_5$
and the nondegeneracy and unfolding conditions for the ordinary and
degenerate flip are invariant under cyclic relabelling.
\end{theorem}

\begin{proof}
The doubled defining sequence has period $q$. Reindexing its $2q$
equations by $k$ therefore gives
\[
    \widehat\Phi^{[k]}(\widehat\rho_k\widehat X,\lambda)
      =\widehat\rho_k\widehat\Phi(\widehat X,\lambda).
\]
Differentiating at the repeated orbit gives
\[
    \widehat L^{[k]}
      =\widehat\rho_k\widehat L\widehat\rho_k^{-1}.
\]
Since $\widehat\rho_k$ is orthogonal, it transports the right and
left kernels, their orthogonal complements and the projection
$\widehat Q$ to the shifted system. Applying the transformed
complementary equation to
$\widehat\rho_k\widehat Z(s,\lambda)$ and using uniqueness gives
\[
    \widehat Z^{[k]}(s,\lambda)
       =\widehat\rho_k\widehat Z(s,\lambda).
\]
Taking the inner product with the shifted left-kernel vector gives
\[
\begin{aligned}
    \widehat\psi^{[k]}(s,\lambda)
    &=\left\langle\widehat\rho_k\widehat W,
      \widehat\Phi^{[k]}\bigl(
          \widehat\rho_k(\widehat X_0(\lambda)
                 +s\widehat V+\widehat Z(s,\lambda)),\lambda\bigr)
       \right\rangle\\
    &=\left\langle\widehat\rho_k\widehat W,
      \widehat\rho_k\widehat\Phi\bigl(
          \widehat X_0(\lambda)
                 +s\widehat V+\widehat Z(s,\lambda),\lambda\bigr)
       \right\rangle
      =\widehat\psi(s,\lambda).
\end{aligned}
\]
Equality of the coefficients and their parameter derivatives now
follows from \eqref{eq:beta35}.
\end{proof}

\section{Examples}
\label{sec:examples}

The examples revisit the scalar fold of \citet{Oliveira2018} and
treat cusps and a degenerate flip in planar and periodically
non-autonomous implicit systems.

\subsection{A period-three fold revisited}

Consider the scalar implicit relation from \citep{Oliveira2018},
\begin{equation}
    F(x,y,\alpha)
    =y-\alpha x\left(1-x+\frac{y^5}{100}\right)=0.
    \label{eq:old-fold-F}
\end{equation}
A period-three fold is located by solving the cyclic orbit equations
and the singularity condition
\begin{equation}
    \Phi(X,\alpha)=0,
    \qquad
    \det L_1=0,
    \label{eq:old-fold-system}
\end{equation}
where $X=(x_0,x_1,x_2)$. The solution is, to eight decimal places,
\begin{equation}
\begin{split}
    x_0^*&=0.16498025,\\
    x_1^*&=0.51812976,\\
    x_2^*&=0.95400016,\\
    \alpha_*&=3.75937510.
\end{split}
    \label{eq:old-fold-orbit}
\end{equation}
The cyclic linearisation $L=L_1=D_X\Phi(X_*,\alpha_*)$ is
\begin{equation}
    L=
    \begin{pmatrix}
    -2.52033362 & 0.99776503 & 0\\
    0 & 0.10660618 & 0.91932890\\
    0.99986715 & 0 & 3.41350920
    \end{pmatrix},
    \label{eq:old-fold-L}
\end{equation}
and the monodromy multiplier is $M=J_2J_1J_0=1$. With the unit
right- and left-kernel vectors
\begin{equation}
\begin{split}
    V&=(0.36596981,\ 0.92443209,\ -0.10719795)^{\mathsf T},\\
    W&=(0.10262415,\ -0.96049580,\ 0.25868146)^{\mathsf T},
\end{split}
    \label{eq:old-fold-VW}
\end{equation}
Propositions~\ref{prop:a2a3} and \ref{prop:param-derivatives} give
\begin{equation}
    \psi_{ss}(0,\alpha_*)=-6.07428006,
    \qquad
    \psi_{\alpha}(0,\alpha_*)=0.21824457.
    \label{eq:old-fold-coeffs}
\end{equation}
The two nonzero derivatives verify the fold conditions of
Theorem~\ref{thm:fold}. Under each of the three cyclic choices of the
initial state, the simultaneous permutation of $X$, $V$ and $W$
preserves these coefficients exactly, as established by
Theorem~\ref{thm:cyclic-invariance}. This recovers the scalar
calculation of \citet{Oliveira2018} within the cyclic formulation.

\subsection{A coupled vector-valued implicit cusp}

Let $(a,b)\in\R^2$ be parameters and $h>0$ the step size of
backward Euler applied to the planar system
\begin{equation}
    \dot u=a+bu+u^3+v+u^5,
    \qquad
    \dot v=-v+u^3.
    \label{eq:vector-cusp-ode}
\end{equation}
With $x=(u,v)$ and $y=(\bar u,\bar v)$ denoting successive states,
the implicit step is
\begin{equation}
F(x,y;a,b)=
\begin{pmatrix}
 \bar u-u-h\bigl(a+b\bar u+\bar u^3+\bar v+\bar u^5\bigr)\\
 \bar v-v+h\bar v-h\bar u^3
\end{pmatrix}=0.
\label{eq:vector-cusp-implicit}
\end{equation}
For a fixed point ($q=1$), the cyclic orbit map becomes
\begin{equation}
    \Phi(u,v;a,b)
    =
    \begin{pmatrix}
      -h(a+bu+u^3+v+u^5)\\
      h(v-u^3)
    \end{pmatrix}.
    \label{eq:vector-cusp-Phi}
\end{equation}
At $(u,v,a,b)=(0,0,0,0)$,
\begin{equation}
    A_0=-I,
    \qquad
    B_0=
    \begin{pmatrix}
       1&-h\\0&1+h
    \end{pmatrix},
    \qquad
    M=J_0=
    \begin{pmatrix}
       1&\dfrac{h}{1+h}\\[1mm]
       0&\dfrac{1}{1+h}
    \end{pmatrix}.
    \label{eq:vector-cusp-linear}
\end{equation}
Thus $+1$ is algebraically simple, and the other multiplier is
$(1+h)^{-1}\in(0,1)$. The cyclic linearisation is
\begin{equation}
    L=A_0+B_0=
    \begin{pmatrix}
       0&-h\\0&h
    \end{pmatrix}.
    \label{eq:vector-cusp-L}
\end{equation}
Choose
\begin{equation}
    V=(1,0)^{\mathsf T},
    \qquad
    W=2^{-1/2}(1,1)^{\mathsf T}.
    \label{eq:vector-cusp-VW}
\end{equation}
Writing $X=sV+Z$ with $Z=(0,z)$, the range equation
$Q\Phi(X;a,b)=0$ reduces to $a+bs+2z+s^5=0$. Hence
\begin{equation}
    z(s;a,b)=-\frac12\bigl(a+bs+s^5\bigr),
    \label{eq:vector-cusp-z}
\end{equation}
and the scalar equation is obtained exactly:
\begin{equation}
    \psi(s;a,b)
    =-\frac{h}{\sqrt2}\bigl(a+bs+2s^3+s^5\bigr).
    \label{eq:vector-cusp-reduced}
\end{equation}
At the origin,
\begin{equation}
    \psi_{ss}(0;0,0)=0,
    \qquad
    \psi_{sss}(0;0,0)=-6\sqrt2\,h\neq0,
    \label{eq:vector-cusp-deg}
\end{equation}
and
\begin{equation}
    \begin{pmatrix}
    \psi_a & \psi_b\\
    \psi_{sa} & \psi_{sb}
    \end{pmatrix}_{(0;0,0)}
    =-\frac{h}{\sqrt2}
    \begin{pmatrix}
    1&0\\0&1
    \end{pmatrix}.
    \label{eq:vector-cusp-rank}
\end{equation}
Theorem~\ref{thm:cusp} therefore gives a cusp of codimension two.
The same polynomial $a+bu+2u^3+u^5$ follows from the equilibrium
equations of \eqref{eq:vector-cusp-ode} by eliminating $v$ through $v=u^3$.
Thus backward Euler carries this equilibrium cusp to a fixed-point
cusp for every $h>0$, with a contracting transverse multiplier.

\paragraph{A nonzero cubic correction.}
For a separate planar example, consider $x=(u,v)$,
$y=(\bar u,\bar v)$ and the implicit relation
\[
G(x,y;a,b)=
\begin{pmatrix}
 \bar u-u+\bar v-\bar u^2+\bar u\bar v-a-b\bar u\\
 2\bar v-v-\bar u^2
\end{pmatrix}=0.
\]
For $q=1$, the cyclic map $\Phi_G(x;a,b)=G(x,x;a,b)$ is
\[
\Phi_G(u,v;a,b)=
\begin{pmatrix}v-u^2+uv-a-bu\\v-u^2\end{pmatrix}.
\]
At the origin, $D_2G$ is invertible and the multipliers are $1$
and $1/2$. The corresponding unit right- and left-kernel vectors are
$V_G=(1,0)^{\mathsf T}$ and $W_G=2^{-1/2}(1,-1)^{\mathsf T}$.
Let $\mathcal R_G$ be the corresponding reduced inverse and
$\psi_G$ the reduced scalar equation. Then
\[
\begin{gathered}
 D_X^2\Phi_G[V_G,V_G]=(-2,-2)^{\mathsf T},\qquad
 \mathcal R_G D_X^2\Phi_G[V_G,V_G]=(0,-2)^{\mathsf T},\\
 D_X^3\Phi_G=0,\qquad
 D_X^2\Phi_G[V_G,(0,-2)^{\mathsf T}]=(-2,0)^{\mathsf T}.
\end{gathered}
\]
Formula~\eqref{eq:psi-sss} therefore gives
$\partial_s^3\psi_G(0;0,0)=3\sqrt2$, entirely from the
complementary correction. Independently, the range equation gives
$v=(2s^2+a+bs)/(2+s)$ and hence the exact reduced equation
\[
 \psi_G(s;a,b)=\frac{\sqrt2\,(s^3-a-bs)}{2+s}.
\]
Eliminating $v$ directly from $\Phi_G=0$ gives
$v=u^2$ and $u^3-a-bu=0$, confirming the calculation.

\subsection{A degenerate flip in a 2-periodic non-autonomous system}

For scalar states $u,v$ and parameters $\lambda=(\mu,\nu)\in\R^2$,
consider the alternating relations
\begin{equation}
\begin{aligned}
 F_0(u,v;\lambda)&=(v-1)+(v-1)^3-u,\\
 F_1(u,v;\lambda)&=v+v^3+(1+\mu)(u-1)+(2+\nu)(u-1)^3,
 \qquad F_{j+2}=F_j.
\end{aligned}
\label{eq:gpd-example-F}
\end{equation}
Both next-state derivatives are positive. The orbit $(0,1)$ has
least period two for all parameters, and its multiplier is
\begin{equation}
 M(\lambda)=-(1+\mu).
 \label{eq:gpd-example-multiplier}
\end{equation}
Thus $p=q=2$, and the doubled system has four scalar states. At
$\lambda=0$, its repeated orbit is $\widehat X_*=(0,1,0,1)$.
Unit right- and left-kernel vectors are
\begin{equation}
 \widehat V=\tfrac12(1,1,-1,-1)^{\mathsf T},\qquad
 \widehat W=\tfrac12(1,-1,-1,1)^{\mathsf T}.
 \label{eq:gpd-example-linear}
\end{equation}
Oddness in the deviations from $\widehat X_*$, exchange symmetry
and uniqueness of the range solution imply that this solution has
the form $(a,1+b,-a,1-b)$. The critical coordinate is $s=a+b$.
Writing $a=s/2+z$ and $b=s/2-z$, the range equation becomes
\begin{equation}
 -2a-a^3+(2+\mu)b+(3+\nu)b^3=0,
 \qquad a=s/2+z,\quad b=s/2-z.
 \label{eq:gpd-example-range}
\end{equation}
Its derivative with respect to $z$ at the origin is $-4$, so it
defines a unique smooth $z=z(s,\lambda)$. The complementary
correction and the exact reduced equation are
\begin{equation}
 \widehat Z=z(1,-1,-1,1)^{\mathsf T},\qquad
 \widehat\psi(s;\lambda)=2(b+b^3-a).
 \label{eq:gpd-example-reduced}
\end{equation}
At the critical parameter, expansion of
\eqref{eq:gpd-example-range} gives
\[
 z(s,0)=\frac{s^3}{16}-\frac{3s^5}{64}+O(s^7),\qquad
 \widehat\psi(s;0)=\frac{3s^5}{32}+O(s^7).
\]
More generally, substitution into
\eqref{eq:gpd-example-reduced} gives
\begin{equation}
 \beta_1(\mu,\nu)=-\frac{2\mu}{4+\mu},\qquad
 \beta_3(\mu,\nu)=\frac{\mu}{4}-\frac{\nu}{8}
                  +O(\|\lambda\|^2),\qquad
 \beta_5(0)=\frac{3}{32}.
 \label{eq:gpd-example-betas}
\end{equation}
Consequently,
\[
 D_{(\mu,\nu)}(\beta_1,\beta_3)(0)=
 \begin{pmatrix}-1/2&0\\1/4&-1/8\end{pmatrix},
 \qquad \det D_{(\mu,\nu)}(\beta_1,\beta_3)(0)=\frac1{16}.
\]
Theorem~\ref{thm:flip-gpd} therefore gives a degenerate flip of
the period-two branch, with bifurcating period-four orbits. Both
defining relations are cubic: the nonzero fifth-order coefficient
comes entirely from the complementary correction.

Changing the starting phase rotates the four states by one place.
The transported kernel vectors are
$\widehat V^{[1]}=\widehat W$ and
$\widehat W^{[1]}=-\widehat V$. The range solution becomes
$(1+b,-a,1-b,a)$; projection onto $\widehat W^{[1]}$ again gives
$2(b+b^3-a)$. Thus the two phases give exactly the same reduced
equation, including its unfolding coefficients.


\subsection{A coupled period-two cusp in a planar implicit system}
\label{sec:alternating-cusp}

We now consider a cusp of period-two orbits in a coupled planar system
with alternating defining relations.  Write
$x_j=(u_j,y_j)\in\R^2$, set $\delta_0=u_0$ and $\delta_1=u_1-1$, and let
$(a,b)\in\R^2$ be the bifurcation parameters.  Define
\begin{equation}
    \mathcal N(z;a,b)=a+bz+2z^3+z^5.
    \label{eq:alternating-nonlinearity}
\end{equation}
The two implicit steps are
\begin{equation}
    F_0(x_0,x_1;a,b)=
    \begin{pmatrix}
       \delta_1-\delta_0-\mathcal N(\delta_1;a,b)-y_1\\
       y_1-\tfrac12 y_0-\tfrac12 \delta_1^3
    \end{pmatrix}=0,
    \label{eq:alternating-F0}
\end{equation}
and
\begin{equation}
    F_1(x_1,x_0;a,b)=
    \begin{pmatrix}
       \delta_0-\delta_1-\mathcal N(\delta_0;a,b)-y_0\\
       y_0-\tfrac13 y_1-\tfrac23 \delta_0^3
    \end{pmatrix}=0.
    \label{eq:alternating-F1}
\end{equation}
The variables are coupled in both directions: the first component
responds to the second coordinate, while the second is forced by the cube of $\delta$.
The relations alternate with period two.  At $(a,b)=(0,0)$ they
admit the orbit
\[
    X_*=\bigl((0,0),(1,0)\bigr),
\]
whose states are distinct.

The second-argument Jacobians at this orbit are both
$\begin{psmallmatrix}1&-1\\0&1\end{psmallmatrix}$, and therefore
invertible.  The derivatives of the forward steps and their
monodromy matrix are
\begin{equation}
\begin{aligned}
    J_0&=\begin{pmatrix}1&\tfrac12\\0&\tfrac12\end{pmatrix},
    & J_1&=\begin{pmatrix}1&\tfrac13\\0&\tfrac13\end{pmatrix},
    & M=J_1J_0&=\begin{pmatrix}1&\tfrac23\\0&\tfrac16\end{pmatrix}.
\end{aligned}
\label{eq:alternating-M}
\end{equation}
Thus $+1$ is algebraically simple and the other multiplier is $1/6$.
In the orbit-space order $(u_0,y_0,u_1,y_1)$, the unit right- and
left-kernel vectors of the cyclic linearisation are
\begin{equation}
    V=\frac{1}{\sqrt2}(1,0,1,0)^{\mathsf T},
    \qquad
    W=\frac{1}{\sqrt{195}}(5,8,5,9)^{\mathsf T}.
    \label{eq:alternating-VW}
\end{equation}
At the critical point, the second state derivatives of $\Phi$
vanish.  Propositions~\ref{prop:a2a3} and
\ref{prop:param-derivatives} give
\begin{equation}
    \psi_{ss}(0;0,0)=0,\qquad
    \psi_{sss}(0;0,0)=-\frac{90}{\sqrt{390}}\ne0,
    \label{eq:alternating-cubic}
\end{equation}
and
\begin{equation}
\begin{pmatrix}
    \psi_a&\psi_b\\
    \psi_{sa}&\psi_{sb}
\end{pmatrix}_{(0;0,0)}
=
\begin{pmatrix}
    -\dfrac{10}{\sqrt{195}}&0\\[1mm]
    0&-\dfrac{10}{\sqrt{390}}
\end{pmatrix}.
\label{eq:alternating-cusp-conditions}
\end{equation}
Both conditions of Theorem~\ref{thm:cusp} hold.  The complementary
correction is nonzero:
$\Phi(X_*+sV;0,0)$ does not lie in the span of $W$ for nonzero $s$
near the origin.

\begin{figure}[htbp]
\centering
\begin{tikzpicture}[>=Latex]
  \draw[->,thin] (-3.50,0) -- (3.55,0) node[right] {$a$};
  \draw[->,thin] (0,-3.20) -- (0,0.56) node[above] {$b$};
  \draw[line width=0.85pt]
    plot[domain=-0.79:0.79,samples=160,variable=\z]
    ({0.73*(6*\z*\z*\z+4*\z*\z*\z*\z*\z)},
     {0.39*(-9*\z*\z-5*\z*\z*\z*\z)});
  \fill (0,0) circle (1.15pt);
\end{tikzpicture}
\caption{Fold locus in the $(a,b)$-parameter plane for the coupled
alternating implicit system.  Its two branches meet at the origin,
where the fold degenerates to a cusp.  Locally, the region between
them contains three period-two solutions, whereas the exterior
contains one.}
\label{fig:alternating-fold-locus}
\end{figure}
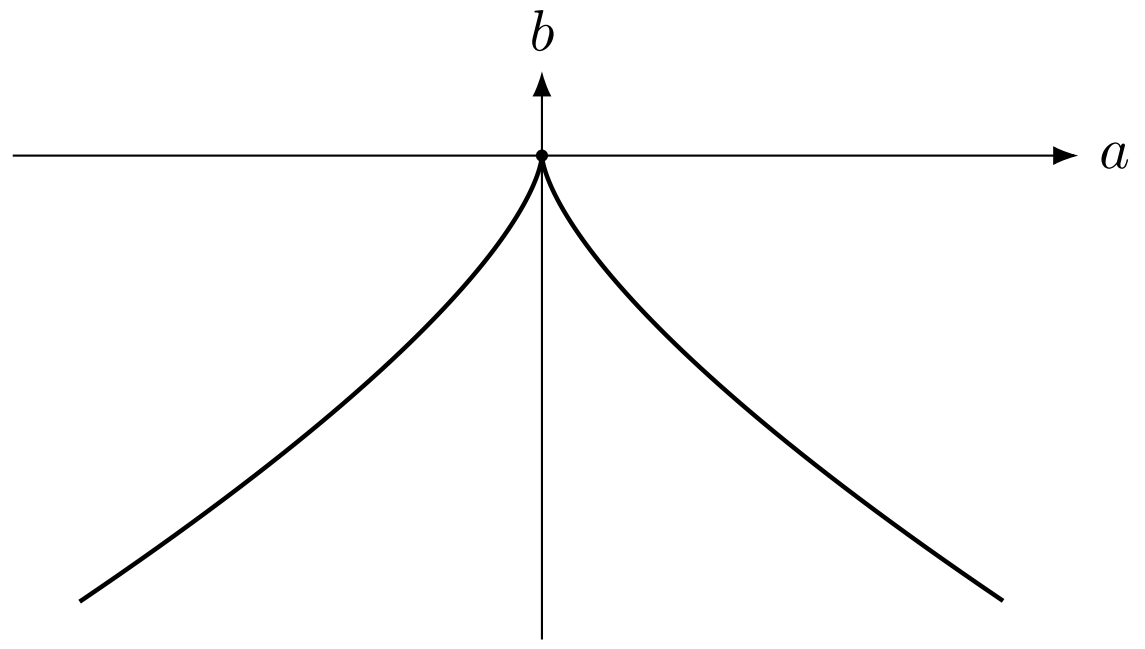

The bifurcation set can also be found exactly.  Put
$\delta_0=t+\eta$ and $\delta_1=t-\eta$.  Use the two $y$-equations and the
difference between the first components of $F_0$ and $F_1$ to
solve for $(\eta,y_0,y_1)$.  At the critical point their Jacobian
with respect to these three variables is
\[
    \begin{pmatrix}
       -4&1&-1\\
       0&-\tfrac12&1\\
       0&1&-\tfrac13
    \end{pmatrix},
    \qquad \det=\frac{10}{3}\ne0.
\]
These three equations have the solution
$\eta=0$, $y_0=y_1=t^3$ for all nearby $(t,a,b)$.
The implicit function theorem makes it the unique local solution.
Substitution into either remaining equation gives
\begin{equation}
    \mathcal C(t;a,b)=a+bt+3t^3+t^5=0.
    \label{eq:alternating-reduced}
\end{equation}
Thus the nearby periodic orbits are precisely
$((t,t^3),(1+t,t^3))$ with $\mathcal C(t;a,b)=0$.
Its solutions correspond to those of the Lyapunov--Schmidt equation,
with $s=\sqrt2\,t$ along the periodic solutions.
The fold locus $\mathcal C=\mathcal C_t=0$ is
\begin{equation}
    a=6t^3+4t^5,\qquad b=-9t^2-5t^4,
    \label{eq:alternating-fold-locus}
\end{equation}
and is shown in Figure~\ref{fig:alternating-fold-locus}.

The change of initial state from $x_0$ to $x_1$ exchanges the
ordered relations $(F_0,F_1)$ and $(F_1,F_0)$.  With the kernel
vectors permuted accordingly, Theorem~\ref{thm:cyclic-invariance}
gives the exact identity $\psi^{[1]}(s;a,b)=\psi(s;a,b)$.

\FloatBarrier

\section{Conclusion}
\label{sec:discussion}

Treating a complete periodic orbit as a cyclic system replaces
differentiation of successive compositions by linear solves with one
reduced operator. The comparison results identify the bifurcation
conditions preserved by regular changes of implicit representation;
cyclic relabelling preserves the reduced equation itself when the
reduction data are transported with the orbit.

A direction for future work concerns bifurcations with
two-dimensional centre dynamics, particularly the discrete
Bogdanov--Takens ($1\!:\!1$ resonance) and Neimark--Sacker
bifurcations. The twisted cyclic operator already identifies the
relevant multipliers, suggesting that the present approach can be
extended to these cases through a two-dimensional reduction.
I leave this question open.

\section*{Declarations}

\paragraph{Funding.}
This work was supported by Funda\c{c}\~ao para a Ci\^encia e a Tecnologia (FCT-PT), through the Center for Mathematical Analysis, Geometry and
Dynamical Systems (CAMGSD), project UID/04459/2025, and by the Recovery and
Resilience Plan and the European Union -- NextGenerationEU through project
UID/PRR/04459/2025.

\paragraph{Competing interests.}
The author declares that he has no competing interests.

\paragraph{Author contributions.}
Single authorship.

\paragraph{Data availability.}
No datasets were generated or analysed in the theoretical development of
this work.

\paragraph{Code availability.}
Wolfram Mathematica was used to check the symbolic and numerical
calculations in the examples. These calculations are elementary and
can be reproduced independently from the equations provided in the
text. The corresponding Mathematica calculations are available from
the author upon reasonable request.

\paragraph{Ethics approval and consent to participate.}
Not applicable.

\paragraph{Consent for publication.}
Not applicable.

\paragraph{Use of generative artificial intelligence.}
During the preparation of the manuscript, the author used OpenAI's ChatGPT
for language editing and the revision of the final draft. All mathematical content,
arguments, calculations and conclusions were written by the author, who
takes full responsibility for the final manuscript.

\section*{ORCID}
\noindent Henrique M. Oliveira:
\url{https://orcid.org/0000-0002-3346-4915}

\bibliographystyle{plainnat}
\bibliography{LS_Implicit_References}

\end{document}